\documentclass[final]{article}
\usepackage{amsmath,amssymb,amsthm,enumerate}
\usepackage{graphicx}
\usepackage[left=1in,top=1in,right=1in,bottom=1in,letterpaper]{geometry}
\usepackage{setspace,url,epsfig}
\usepackage[numbers]{natbib}
\usepackage{algorithm,algorithmic}
\usepackage{subfigure}
\newcommand{\sign}{\text{sign}}
\newcommand{\be}{\begin{equation}}
\newcommand{\ee}{\end{equation}}
\newcommand{\benn}{\begin{equation*}}
\newcommand{\eenn}{\end{equation*}}

\newtheorem{proposition}{Proposition}

\begin{document}

\title{A nonlinear PDE-based method for sparse deconvolution
\thanks{Research supported by ONR grants: N00014-08-1-1119 and N00014-07-0810 and the Department of Defense}}
\author{
 Yu Mao
 \thanks{Department of Mathematics, UCLA, Los Angeles, CA 90095 (\texttt{ymao29@math.ucla.edu})}
 \and Bin Dong
 \thanks{Department of Mathematics, UCSD, 9500 Gilman Drive, La Jolla,
CA, 92093-0112  , (\texttt{b1dong@math.ucsd.edu})}
 \and Stanley Osher
 \thanks{Department of Mathematics, UCLA, Los Angeles, CA 90095 (\texttt{sjo@math.ucla.edu})}
}
\date{February 10, 2010}

\maketitle

\begin{abstract}
In this paper, we introduce a new nonlinear evolution partial
differential equation for sparse deconvolution problems. The
proposed PDE has the form of continuity equation that arises in
various research areas, e.g. fluid dynamics and optimal
transportation, and thus has some interesting physical and
geometric interpretations. The underlying optimization model that
we consider is the standard $\ell_1$ minimization with linear
equality constraints, i.e. $\min_u\{\|u\|_1 : Au=f\}$ with $A$
being an under-sampled convolution operator. We show that our PDE
preserves the $\ell_1$ norm while lowering the residual
$\|Au-f\|_2$. More importantly the solution of the PDE becomes
sparser asymptotically, which is illustrated numerically.
Therefore, it can be treated as a natural and helpful plug-in to
some algorithms for $\ell_1$ minimization problems, e.g. Bregman
iterative methods introduced for sparse reconstruction problems in
\cite{Yin:2008p611}. Numerical experiments show great improvements
in terms of both convergence speed and reconstruction quality.
\end{abstract}

\section{Introduction}\label{introduction}

The sparse deconvolution problem is to estimate an unknown sparse
signal which has been convolved with some known kernel function
and corrupted by noise. There are many applications of sparse
deconvolution, e.g. in seismology, nondestructive ultrasonic
evaluation, image restoration, or intracardiac electrograms (see
e.g.
\cite{mendel1990maximum,o1994recovery,olofsson2007sparse,olofsson1999maximum,olofsson2004semi,galatsanos1991least,ellis1996deconvolution}).

There has been an intensive development of sparse deconvolution
techniques and algorithms. Due to their simplicity, least squares methods have been used to deconvolve sparse seismic signals
\cite{berkhout1977least,LawsonLSQ}. The major drawbacks of least squares
methods are their lack of robustness and the ill-poseness of the
problem when the number of unknowns is less than or equal to the
number of equations. Many complementary methods have been
presented, such as Tikhonov regularization \cite{tikhonov1977}
and the total least squares method \cite{golub1980analysis,van1991total}. Another
class of complementary methods to the classical least squares method is the
$\ell_1$-penalized models
\cite{taylor1979deconvolution,debeye1990lp}, where $\ell_1$
norm of the unknown signal is used as an additional regularization
term, allowing us to control the sparseness of the solution.

In the recent burst of research in compressive sensing (CS)
\cite{Candes:2006p1807,Candes:2006p625,Donoho:2006p668}, the
following $\ell_1$ minimization problem (basis pursuit) is widely
used (see e.g. \cite{CDS_BP,Yin:2008p611})
\be\label{L1:Minimization} \min_{u \in R^n}\{\|u\|_1 : Au=f \},
\ee where $\|u\|_1:=\sum_i|u_i|$ denotes the $\ell_1$ norm of $u$.
In our paper, we shall adopt the above optimization model with $A$
being a convolution matrix with some known continuous differentiable kernel function (e.g. the impulse response of the system or wavelet basis). When the kernel
function is everywhere nonnegative and the sparse signal we want
to recover is also known to be nonnegative, we can solve the
constrained least squares problem $\min_u\{\|Au-f\|_2^2: u\ge0\}$ as suggested
by \cite{Bruckstein:2008p2703}. However, this model does not work for
general kernel functions and signals with negative entries. Hence
in this paper, we will stick to the above $\ell_1$ minimization
model \eqref{L1:Minimization}, although we believe our PDE can be also added to many  nonnegative least squares approaches.

To solve the constrained optimization \eqref{L1:Minimization},
standard second-order methods such as interior-point methods were
introduced in \cite{candes20071,koh-solver,lustig2007sparse}.
However, these methods, although accurate, become inefficient for
large scale problems. To overcome this, large scale optimization
methods have recently been developed to solve
\eqref{L1:Minimization} or its siblings
\cite{Yin:2008p611,DO,OMDY,COS1,COS2,zhang2009bregmanized,becker2009nesta,HYZ,figueiredo2007gradient,van2008probing}.
These methods are very popular in CS, and all of them, especially
Bregman iterations introduced to CS in \cite{Yin:2008p611}, have outstanding
performance in solving CS problems, e.g. when the matrix $A$ is
taken to be a random submatrix of some orthonormal matrix, such as
Fourier or DCT matrix.

For deconvolution problems however, the above mentioned methods
are usually less efficient than solving CS problems. The major
reason is that each column vector of a convolution matrix $A$ is
highly coherent to the nearby columns, especially the ones that
are right next to it. This makes the signal much harder to
deconvolve if it has two spikes that are very close to each other.
Therefore, in contrast to CS problems, the spatial information of
$u$ (i.e. the variable ``$x$" of $u(x)$) is no longer irrelevant
for sparse deconvolution problems, which motivates us to introduce what we call
a spatial motion PDE into optimization algorithms to solve problem
\eqref{L1:Minimization}.

Intuitively speaking, this PDE should spatially move spikes around
such that the constraint $Au=f$ is closer to being satisfied.
However, the $\ell_1$ norm of the solution should not be altered
during the process. We will prove in Section \ref{theory} that our
motion PDE indeed satisfies these properties. Furthermore, the PDE
alone tends to sparsify the solution, which, together with its
other properties, makes it an effective and powerful add-on to
some of the optimization algorithms for CS problems mentioned
above. In this paper, we shall combine the motion PDE with Bregman
iteration \cite{Yin:2008p611} (with the inner unconstrained
problem \eqref{Subproblem:Bregman} solved by FPC \cite{HYZ}) to
enhance its performance. We note that it is also possible to
combine the motion PDE with other iterative thresholding based
algorithms (e.g. the linearized Bregman iterative method in
\cite{OMDY} or the Bregmanized operator splitting method in
\cite{zhang2009bregmanized}).

The rest of the paper is organized as follows. We will first
recall some basics of Bregman distance and Bregman iterations in
Section \ref{background}. Then in Section \ref{theory}, we
introduce our motion PDE and explore its physical and geometrical
meaning. In Section \ref{numerical}, we combine our PDE with
Bregman iteration and propose our new algorithm. Numerical
experiments are also conducted. Finally, we draw our conclusions
in Section \ref{conclusion}.

\section{Background}\label{background}

In this section, we briefly recall the concept of Bregman distance
\cite{Br} and Bregman iterations \cite{Yin:2008p611}.

\subsection{Bregman Distance}

The Bregman distance \cite{Br}, based on the convex function $J$,
between points $u$ and $v$, is defined by
\begin{equation*}
D_J^p(u,v) = J(u) - J(v) - \langle p,u-v\rangle
\end{equation*}
where $p \in  \partial J(v)$ is an element in the subgradient of
$J$ at the point $v$. In general $D_J^p(u,v)\not= D_J^p(v,u)$ and
the triangle inequality is not satisfied, so $D_J^p(u,v)$ is not a
distance in the usual sense. However it does measure the closeness
between $u$ and $v$ in the sense that $D_J^p(u,v) \geq 0$ and
$D_J^p(u,v) \geq D_J^p(w,v)$ for all points $w$ on the line
segment connecting $u$ and $v$. Moreover, if $J$ is convex,
$D_J^p(u,v) \geq 0$, if $J$ is strictly convex $D_J^p(u,v) > 0$
for $u \not= v$ and if $J$ is strongly convex, then there exists a
constant $a > 0$ such that
\[
D_J^p(u,v) \geq a\|u-v\|_2^2.
\]

\subsection{Bregman Iterations}

To solve \eqref{L1:Minimization} Bregman iteration was proposed in
\cite{Yin:2008p611}. Given $u^0 = p^0 = 0$, we define: \be
\label{2.2}
\begin{array}{ll}
u^{k+1} &= \arg\min_{u\in R^{n}} \left\{\mu\|u\|_1-\mu\|u^k\|_1 -
\langle u-u^k,p^k\rangle + \frac{1}{2}
\|Au-f\|_2^2\right\}  \\
p^{k+1} &= p^k - A^T (Au^{k+1}-f).
\end{array}
\ee This can be written as, for $J(u)=\mu\|u\|_1$,
\begin{displaymath}
u^{k+1} = \arg\min_{u\in R^{2}} \left\{D_{J}^{p^{k}} (u,u^k) +
\frac{1}{2} \|Au-f\|_2^2\right\}. \nonumber
\end{displaymath}

As shown in \cite{Yin:2008p611}, see also \cite{OBGXY,CHF}, the
Bregman iteration (\ref{2.2}) can be written as, for $f^0 = 0, u^0
= 0$:
\begin{align}
u^{k+1} &= \arg\min_{u\in R^{n}} \left\{\mu\|u\|_1 + \frac{1}{2}
\|Au-f^{k}\|_2^2\right\} \label{Subproblem:Bregman} \\
f^{k+1} &= f^k + f-Au^{k+1}. \nonumber
\end{align}
At each iteration, we will solve the problem
\eqref{Subproblem:Bregman} via the fixed point continuation (FPC)
method proposed in \cite{HYZ}. Other ways of solving the
unconstrained problem \eqref{Subproblem:Bregman} can be found in
e.g.
\cite{figueiredo2007gradient,kim2007method,donoho2006sparse,nesterov2007gradient,van2008probing,figueiredo2003algorithm}
and the references therein. Our improvement will work with most of
these solvers. Now altogether, we have the following FPC Bregman
iterations solving problem \eqref{L1:Minimization}:
\begin{eqnarray}
u^{k+1,N} &\longleftarrow& \left\{\begin{array}{ll} u^{k,l+\frac{1}{2}} & \leftarrow u^{k,l}-\delta A^\top(Au^{k,l}-f^k)\\
u^{k,l+1} & \leftarrow \text{shrink}(u^{k,l+\frac{1}{2}},\mu) \end{array}\right.\quad \text{(FPC)} \nonumber\\
f^{k+1} &=& f^k + f-Au^{k+1,N}. \label{FPC+Bregman}
\end{eqnarray}
Here ``shrink" is the soft thresholding function \cite{Donoho}
defined as
\begin{align*} \text{shrink}  (x,\mu): =
\begin{cases} x-\mu, & \text{if} \ x > \mu \\
0, & \text{if} \ -\mu \leq x \leq \mu \\
x + \mu, & \text{if} \ x < -\mu.
\end{cases}
\end{align*}
Notice that, if we only perform 1 iteration for the FPC in
\eqref{FPC+Bregman}, i.e. $N=1$, this gives us another algorithm,
called Bregmanized operator splitting, that also solves problem
\eqref{L1:Minimization} \cite{zhang2009bregmanized}.

\section{A Spatial Motion PDE}\label{theory}

In this section, we introduce and analyze a spatial motion PDE.
The main purpose of this PDE is to spatially move and properly
combine spikes in a sparse signal $u$, such that the residual
$\|Au-f\|_2$ is reduced while $\|u\|_1$ is preserved. We will also
present some interesting physical and geometric interpretations of
the motion PDE. Throughout this section, $u$ is understood as a
function defined on the image domain $\Omega$, while $A$ is a
composition of an convolution operator and a spatial restriction
operator with respect  to a sampling region $S$, and $A^\top$ is
the conjugate operator of $A$. In other words, \be
(Au)(x)=\chi_S(x)\int_\Omega K(x-y)u(y)dy \ee where $K$ is a differentiable convolution kernel and $\chi_S$ is the
characteristic function of $S$. The other operators in this
section should be understood as functional operators as well.

\subsection{Motivation}\label{motivation}

To analyze the desired spatial motion of the signal, we start from
a characteristic point of view, i.e. observing the moving trajectory
of a point mass and describing its motion. Consider a
simple sparse signal that consists of an isolated spike at
location $c$: \benn\label{prototype1} u=\delta(x-c). \eenn We want
to find the proper spatial motion that decreases the energy
$H(u)=\frac{1}{2}\|Au-f\|_2^2$. To do this, we take the partial
derivative of the energy with respect to the spatial variable $c$
and obtain \benn\label{prototype2} \frac{\partial}{\partial
c}H(u)=-\nabla_x(H'(u))\big|_{x=c}=-\nabla_x(A^\top(Au-f))\big|_{x=c}.
\eenn The above identity means that if we want to minimize $H(u)$,
we can spatially move the spike along the direction
$\nabla_x(A^\top(Au-f))\big|_{x=c}$. In other words, the
trajectory of the moving spike that minimizes $H(u)$ can be
described as $u(x,t)=\delta(x-c(t))$ where $c(t)$ satisfies
\benn\label{prototype3}
c'(t)=-\nabla_x(A^\top(Au-f))\big|_{x=c(t)}. \eenn Similarly, if
the initial signal is a negative spike, i.e. $u=-\delta(x-c)$,
then its correct direction of spatial movement should be
$\frac{\partial}{\partial
c}H(u)=\nabla_x(A^\top(Au-f))\big|_{x=c}$.

Based on this observation, we can formulate a
spatial motion PDE as follows: \be\label{motionPDE1}
u_t+\nabla_x\cdot[u a(u)]=0 \ee where $a(u)$ is the velocity field
of the spatial transport. As elaborated by the examples above, the
spatial velocity field $a(u)$ should be defined as \benn
a(u)=-\nabla_x(A^\top(Au-f))\cdot \sign(u) \eenn Thus
\eqref{motionPDE1} can be equivalently formulated as
\be\label{motionPDE} u_t=\nabla_x\cdot[|u|\nabla_x(A^\top(Au-f))].
\ee For definiteness we attach the initial and boundary conditions
to \eqref{motionPDE} and obtain the following Cauchy problem
\be\label{CauchyPDE}
\begin{cases}u_t=\nabla_x\cdot[|u|\nabla_x(A^\top(Au-f))],\quad (x,t)\in\Omega\times (0,T]\\
u(\cdot,0)=u_0,u|_{\partial\Omega}=0
\end{cases}
\ee where $u_0$ is a continuous initial function.
\eqref{CauchyPDE} is the core spatial motion PDE we will discuss
throughout this paper. A locally integrable function $u$ is called
a weak solution of \eqref{CauchyPDE} if $u$ satisfies the boundary
condition and for any $C^1$ test function $\eta$ the following
identities hold: \be \int_{\Omega\times(0,T]} [\eta_t\cdot u -
|u|\nabla_x\eta \cdot\nabla_x (A^\top(Au-f))]dxdt=0, \ee \be
\lim_{t\to 0^+}\int_\Omega \eta(x,t) u(x,t)dx=\int_\Omega
\eta(x,0)u_0(x)dx. \ee

\subsection{Physical and Geometric Interpretations of PDE \eqref{motionPDE}}\label{related}

The equation \eqref{motionPDE} has a nice physical interpretation,
and is closely related to optimal transportation and nonlinear
dissipative processes. To see this, we assume $u>0$ for now and
regard $u$ as the mass density of a certain fluid. The mass
conservation can be expressed as the following continuity equation
for $u$: \be\label{continuity} u_t+\nabla_x\cdot[u \vec w]=0 \ee
where $\vec w$ is the velocity field that describes the spatial
motion of the fluid. Darcy's law \cite{Vazquez:2007p3731} connects
the velocity field $\vec w$ with the pressure field $p$ by
\be\label{Darcy} \vec w=-\nabla_x p, \ee and thermodynamics tells
us the pressure $p$ is determined by the potential $E(u)$:
\be\label{thermo} p=\partial_u E(u). \ee

For an ideal gas in a homogeneous porous medium, the potential
$E(u)$ is given by the free energy $\int Cu^m$. Thus equations
\eqref{continuity}, \eqref{Darcy} and \eqref{thermo} lead to the
classical porous medium equation \cite{Vazquez:2007p3731}. If we
replace the potential $E(u)$ by our residual $H(u)$, the
combination of \eqref{continuity}, \eqref{Darcy} and
\eqref{thermo} turns out to be our spatial motion equation
\eqref{motionPDE}.

Equation \eqref{motionPDE} can also be understood
as a \textit{gradient flow} of the potential $H(u)$ under the
Wasserstein metric. In \cite{Otto:2001p4508,Otto:2006p4844}, Otto
\emph{et al.} showed that the porous medium equation is a gradient
flow under the Wasserstein metric for some given potential. Our
following interpretation of \eqref{motionPDE} will be in
the same spirit as in \cite{Otto:2001p4508,Otto:2006p4844}.

Wasserstein distance is widely used in optimal transportation
problems, whose fundamental goal is to find the most efficient
plan to transport a density function $u_0$ to another density
function $u_1$ in a mass preserving fashion. Monge's optimal
transportation problem is to find the solution to \be
\inf_{M\#u_0=u_1} \int u_0(x)|x-M(x)|^2dx \ee and the
corresponding optimal value is called Wasserstein distance of
order 2 (see \cite{Villani:2003p4846, Villani:2003p4509} for the
definition of ``$\#$" and more other details). The idea of optimal
transportation has been applied to various problems in image and
signal processing, such as image classification
\cite{Rubner:1998p4851}, registration \cite{Haker:2004p4845} and
segmentation \cite{Ni:2009p4733}.

As pointed out by Otto \cite{Otto:2001p4508}, the Wasserstein
distance of order 2 can be understood as a geodesic distance under
a certain Riemannian structure on the set of density functions.
The metric tensor $g$ on the tangent space at $u$ is formally
defined by \be\label{metric1} \left\langle
s_1,s_2\right\rangle_g=\int u\nabla_x \phi_1\cdot\nabla_x \phi_2
\ee where $s_1$ and $s_2$ are any two infinitesimal variations of
$u$, and $\phi_i, i=1,2,$ are solutions of \be\label{metric2}
s_{i}+\nabla_x\cdot(u\nabla_x\phi_{i})=0. \ee Given the metric
tensor $g$, the set of all density functions forms a Riemannian
manifold $(\mathcal{M},g)$.

Now we can formally interpret the motion PDE \eqref{motionPDE} as
a gradient flow of $H(u)=\frac{1}{2}\|Au-f\|_2^2$ on
$(\mathcal{M},g)$. Under the Riemannian structure, the
corresponding gradient of energy $H(u)$ w.r.t. $g$, denoted by
$\text{grad}_g H(u)$, is defined by the following identity \be
\left\langle \text{grad}_g H(u), \vec
v\right\rangle_g=\partial_{\vec v} H(u) \ee for all vector field
$\vec v$ on $(\mathcal{M},g)$, where $\partial_{\vec v}$ denotes
the directional derivative along $\vec v$. On the other hand, from
\eqref{metric1} and \eqref{metric2} we can see that $\left\langle
s_1,s_2\right\rangle_g=\int \phi_1 s_2$ where $\phi_1$ solves
$s_1+\nabla_x\cdot(u\nabla_x\phi_1)=0$. Therefore, since $u$
solves $u_t+\nabla_x\cdot(u\nabla_x (-\partial_u H(u)))=0$, we
have \be \left\langle u_t, \vec v\right\rangle_g=\int -\partial_u
H(u) \cdot \vec v=-\partial_{\vec v} H(u),\,\forall \vec v \ee
which indicates that $u_t=-\text{grad}_g H(u)$. This elegant
interpretation leads to a substantial understanding of our motion
PDE \eqref{motionPDE}: it gives the most natural way (in the sense
of optimal transport) to spatially move the signal such that the
residual $H(u)$ is decreasing.

\subsection{Properties}\label{properties}

Now we come back to the general case where $u$ is not assumed to
be positive. We have the following fairly well known properties of
equations that resemble \eqref{CauchyPDE} that explain, to some
extent, why this PDE can be utilized to improve the spatial
reconstruction.

\begin{proposition}
If $u(x,t)$, $(x,t)\in\Omega\times (0,T]$  is a weak solution of  \eqref{CauchyPDE} then
$H(u(x,t))=\frac{1}{2}\|Au-f\|_2^2$ is non-increasing over time.
\end{proposition}
\begin{proof}
We need to show that $\dfrac{d}{dt}H(u(x,t))\leq0$. Indeed,
\begin{align*}
\dfrac{d}{dt}H(u(x,t))&=\int H'(u)\cdot u_t dx\\
&=-\int \nabla_x(A^\top (Au-f))\cdot[|u|\nabla_x(A^\top(Au-f))]dx\\
&=-\int |u|[\nabla_x(A^\top(Au-f))]^2dx\leq 0
\end{align*}
\end{proof}

\begin{proposition}
If $u(x,t)$, $(x,t)\in\Omega\times (0,T]$ is continuous weak solution of \eqref{CauchyPDE} then $\|u\|_1$ is
conserved over time.
\end{proposition}
\begin{proof}
For simplicity, we will prove this in the 1D case, while the proof for
higher dimensions is similar. For a fixed time $t$, suppose
$I=[a,b]$ is one of the connected component of $\{u\neq 0\}$, i.e.
$u$ has constant sign within $I$, then we have $u(a,t)=u(b,t)=0$,
so
\begin{align*}
\dfrac{d}{dt}\int_a^b|u|dx&=\int_a^b u_t \sign(u)dx\\
&=\sign(u) \cdot\int_a^b \nabla_x\cdot[|u|\nabla_x(A^\top(Au-f))]\\
&=\sign(u)\cdot [|u|\nabla_x(A^\top(Au-f))]\big|_{x=a}^{x=b}=0
\end{align*}
Therefore overall we have $\frac{d}{dt}\|u\|_1=0$.
\end{proof}

The above two propositions show that, although the total mass of
$u$ is preserved over time, the spatial mass distribution is
adjusted such that the equation $Au=f$ is better satisfied.

\begin{figure}[h]
\centering
    \includegraphics[width=1.0in]{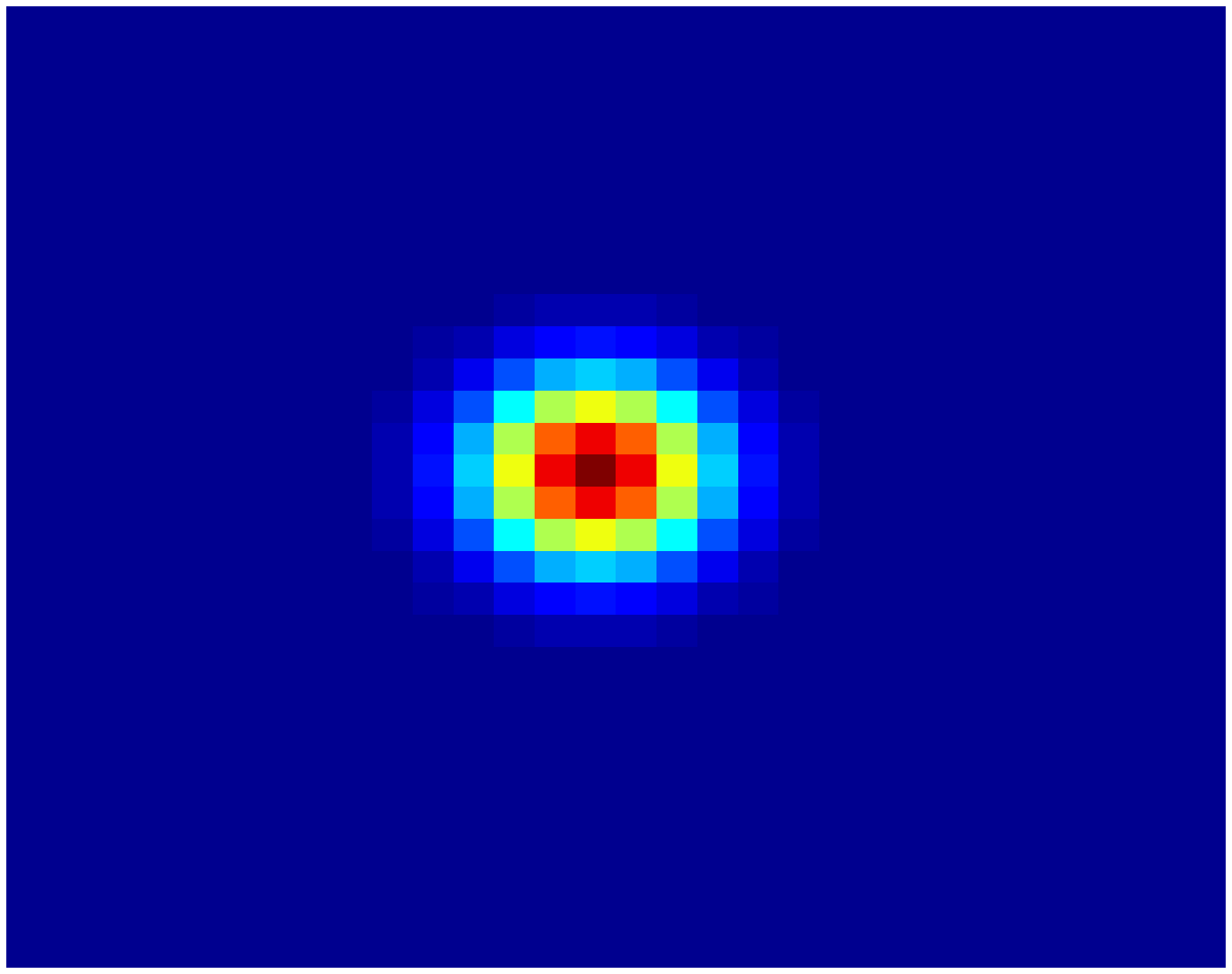}\quad
    \includegraphics[width=1.0in]{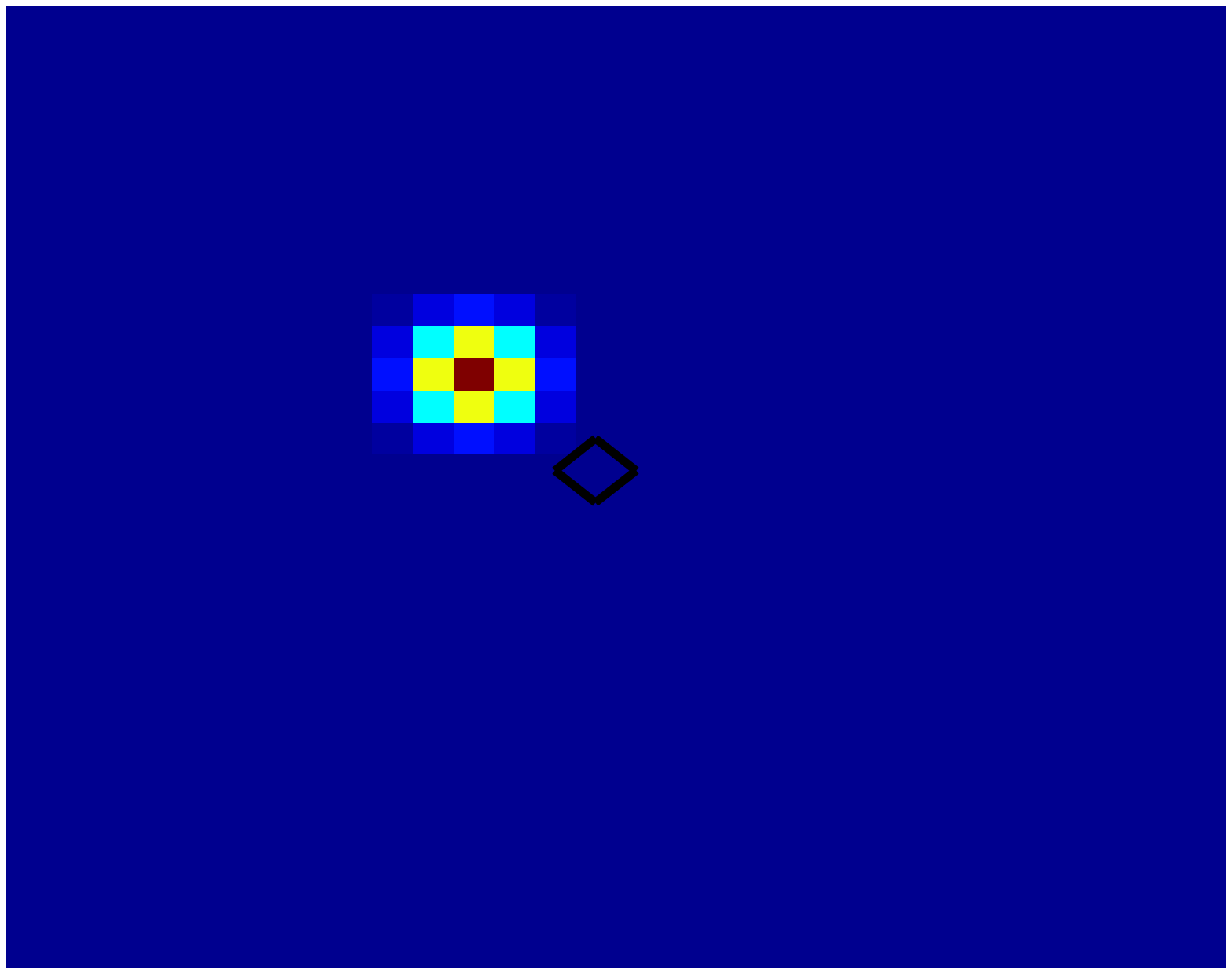}
    \includegraphics[width=1.0in]{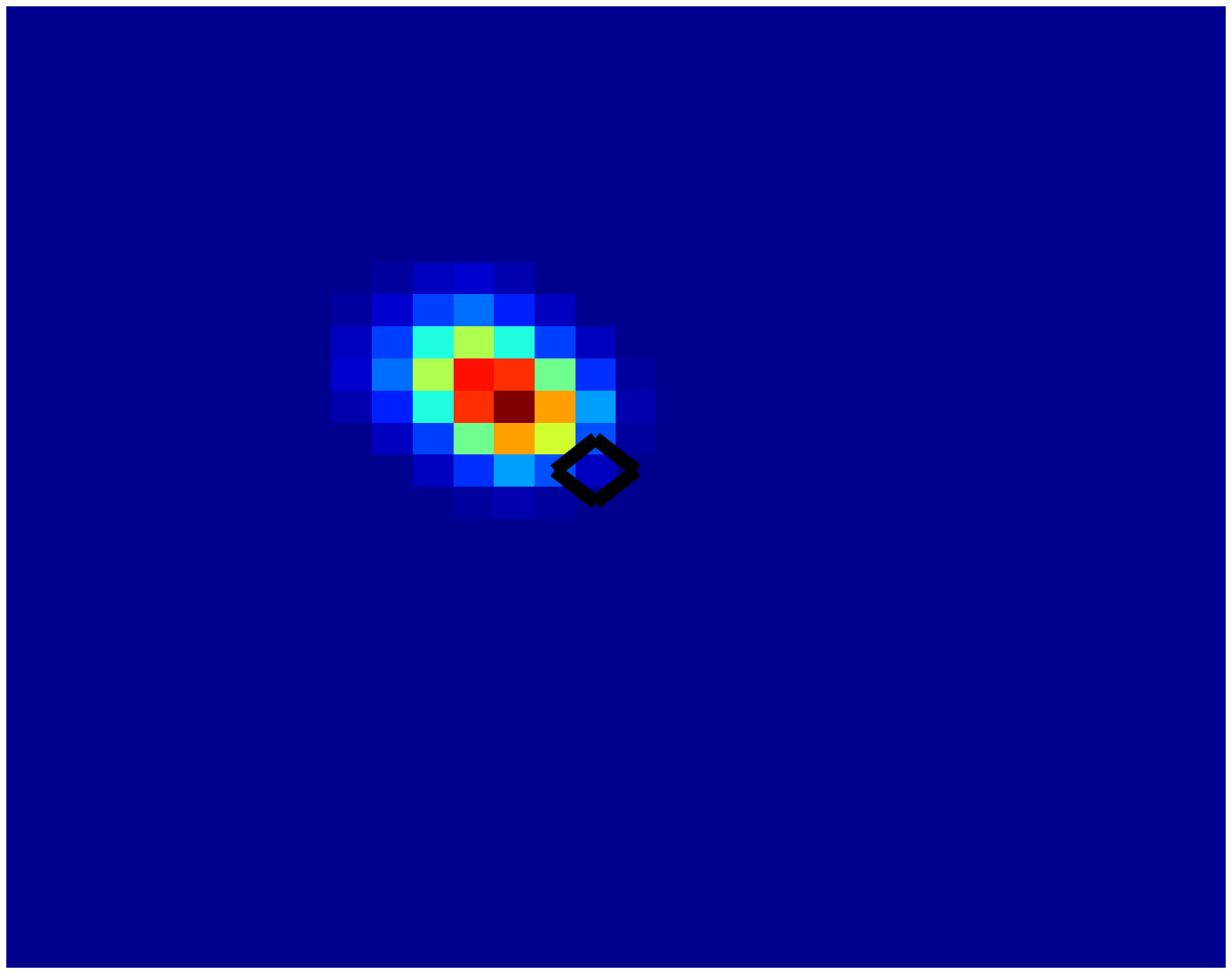}
    \includegraphics[width=1.0in]{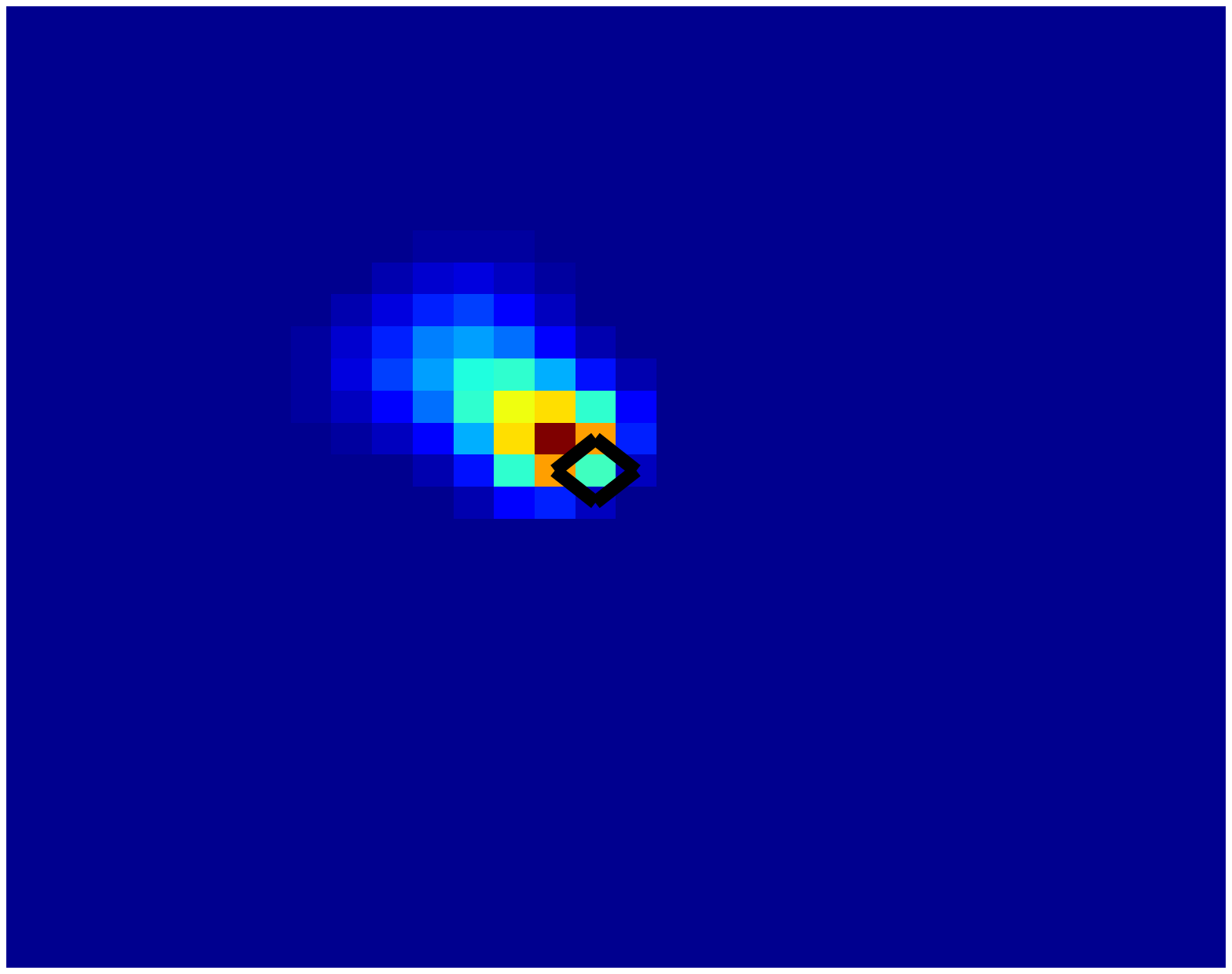}
    \includegraphics[width=1.0in]{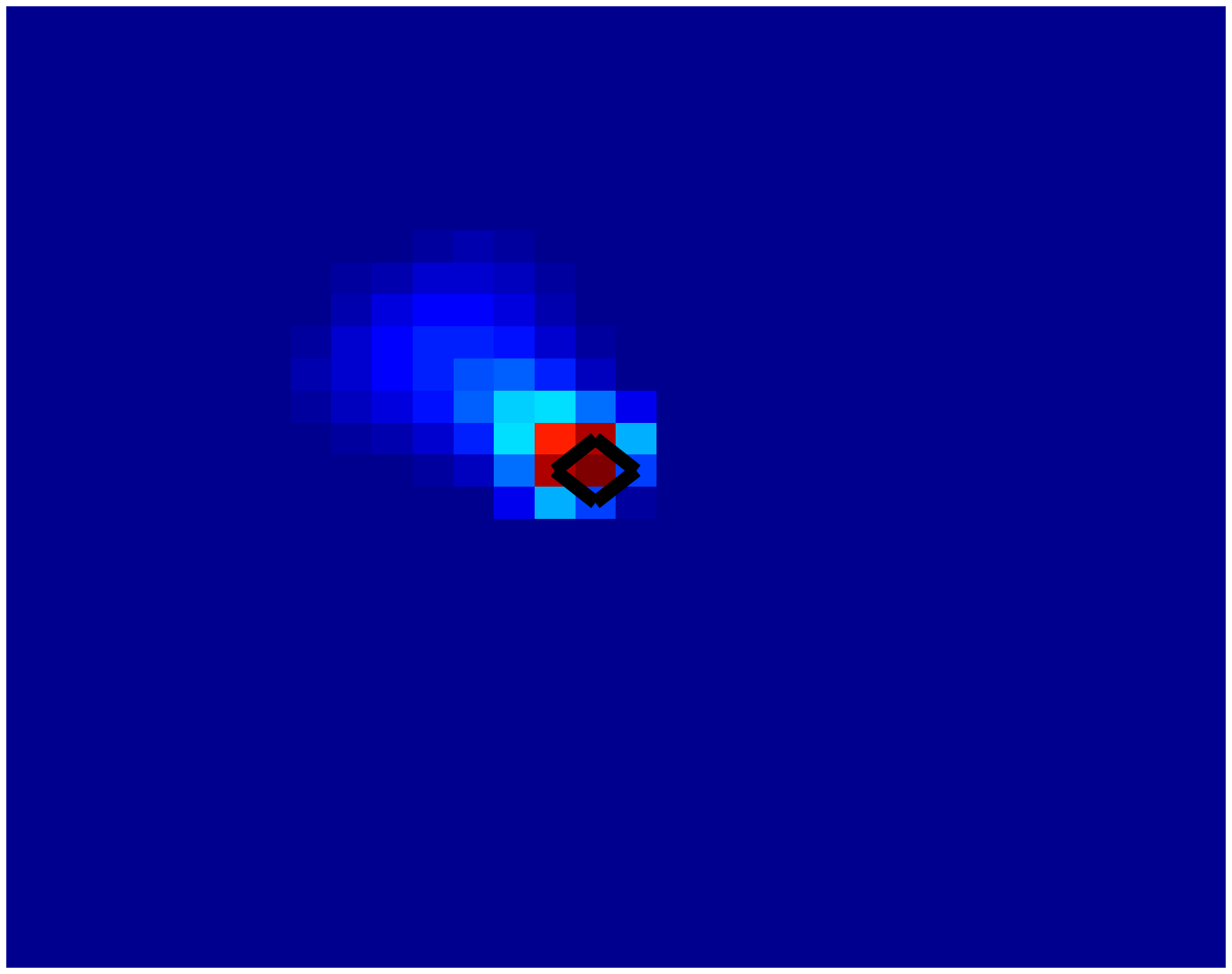}
    \includegraphics[width=1.0in]{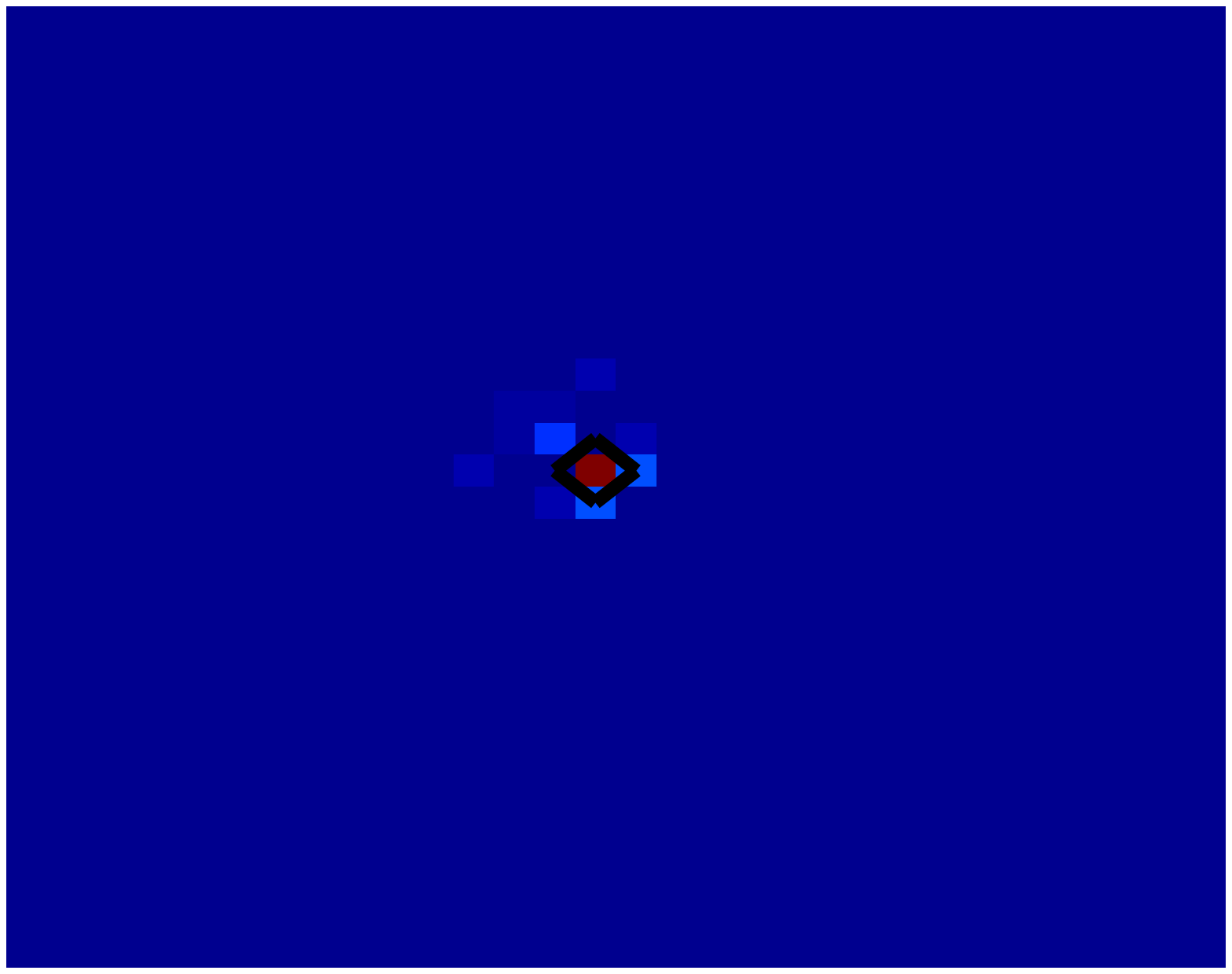}\\
    \includegraphics[width=1.0in]{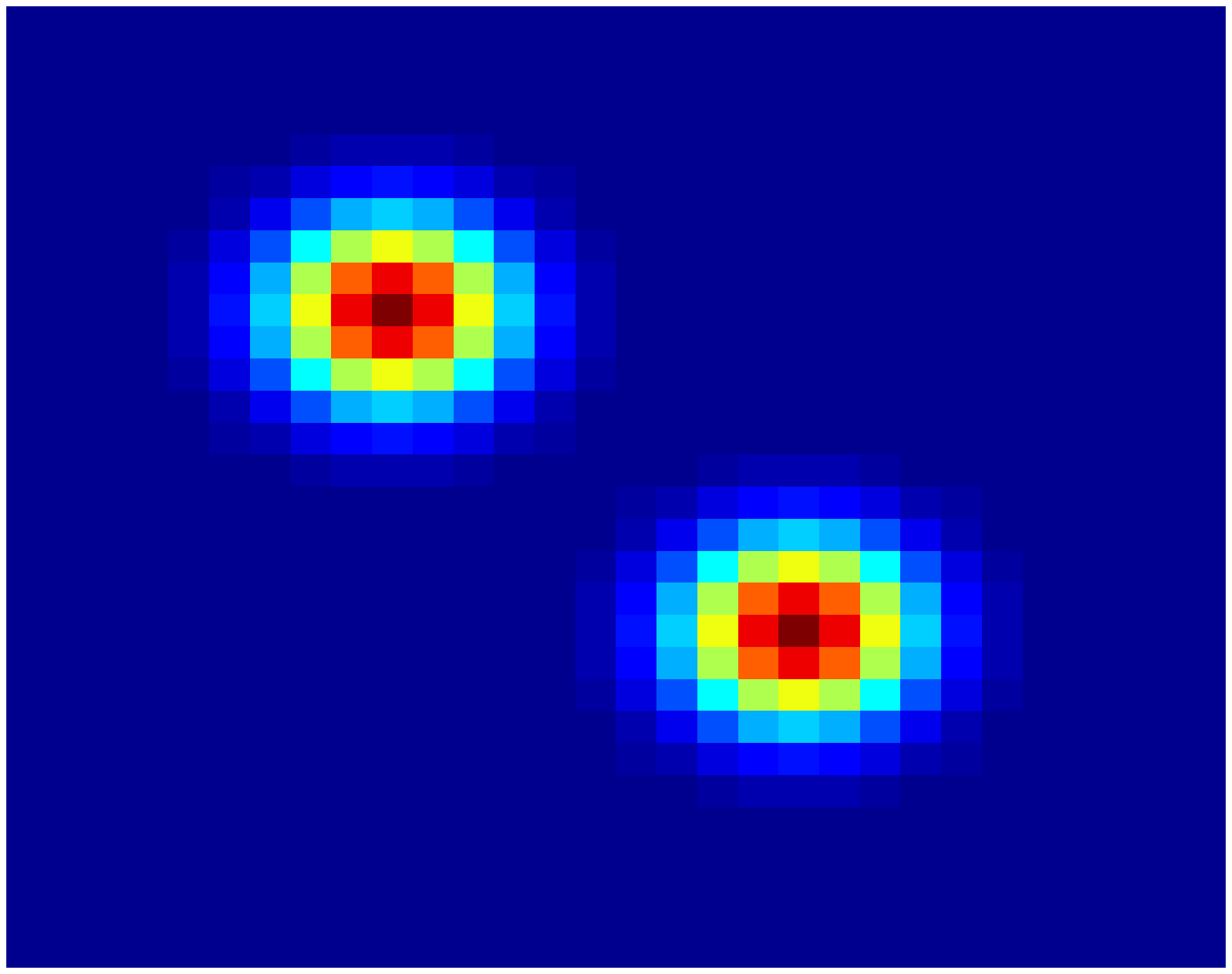}\quad
    \includegraphics[width=1.0in]{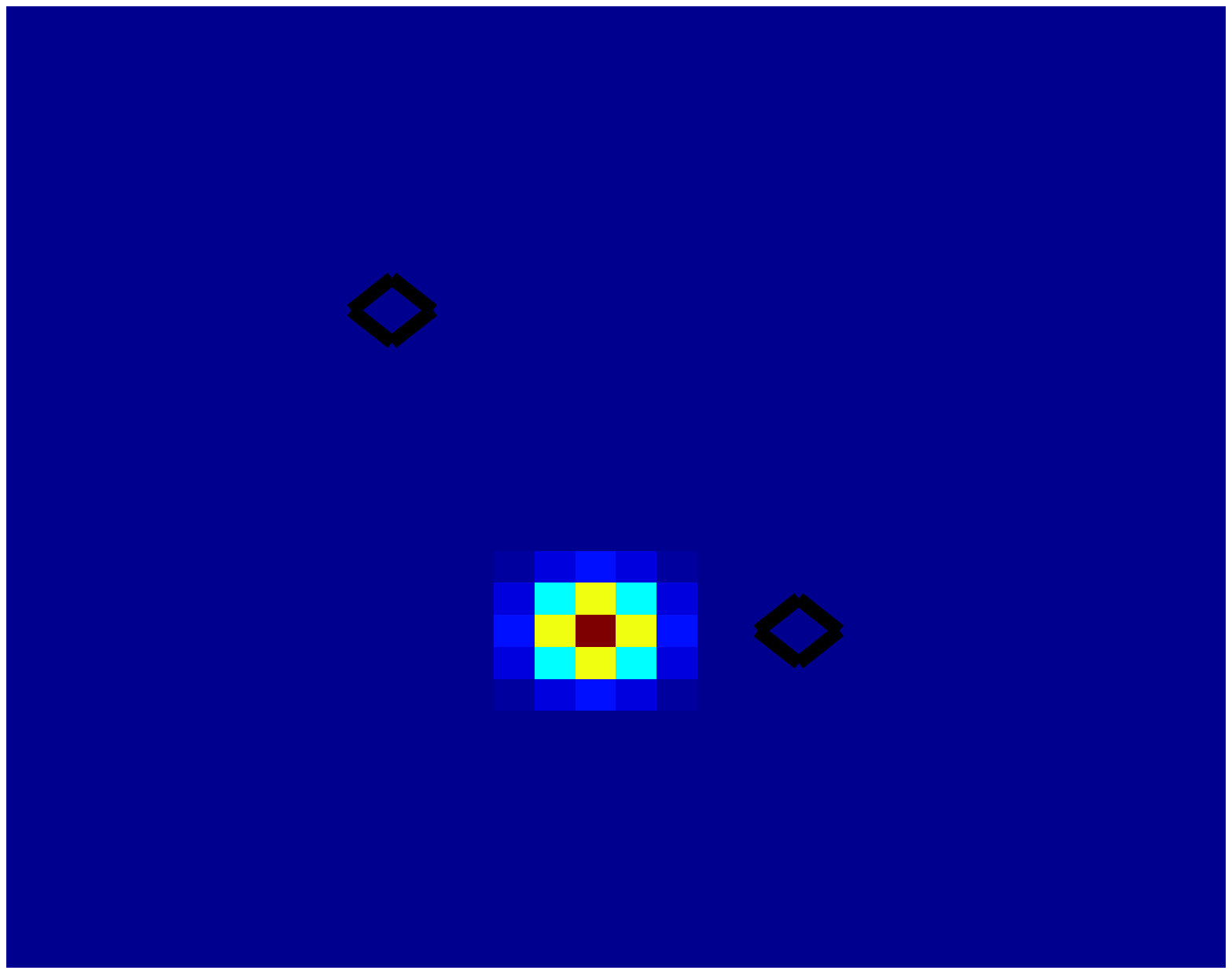}
    \includegraphics[width=1.0in]{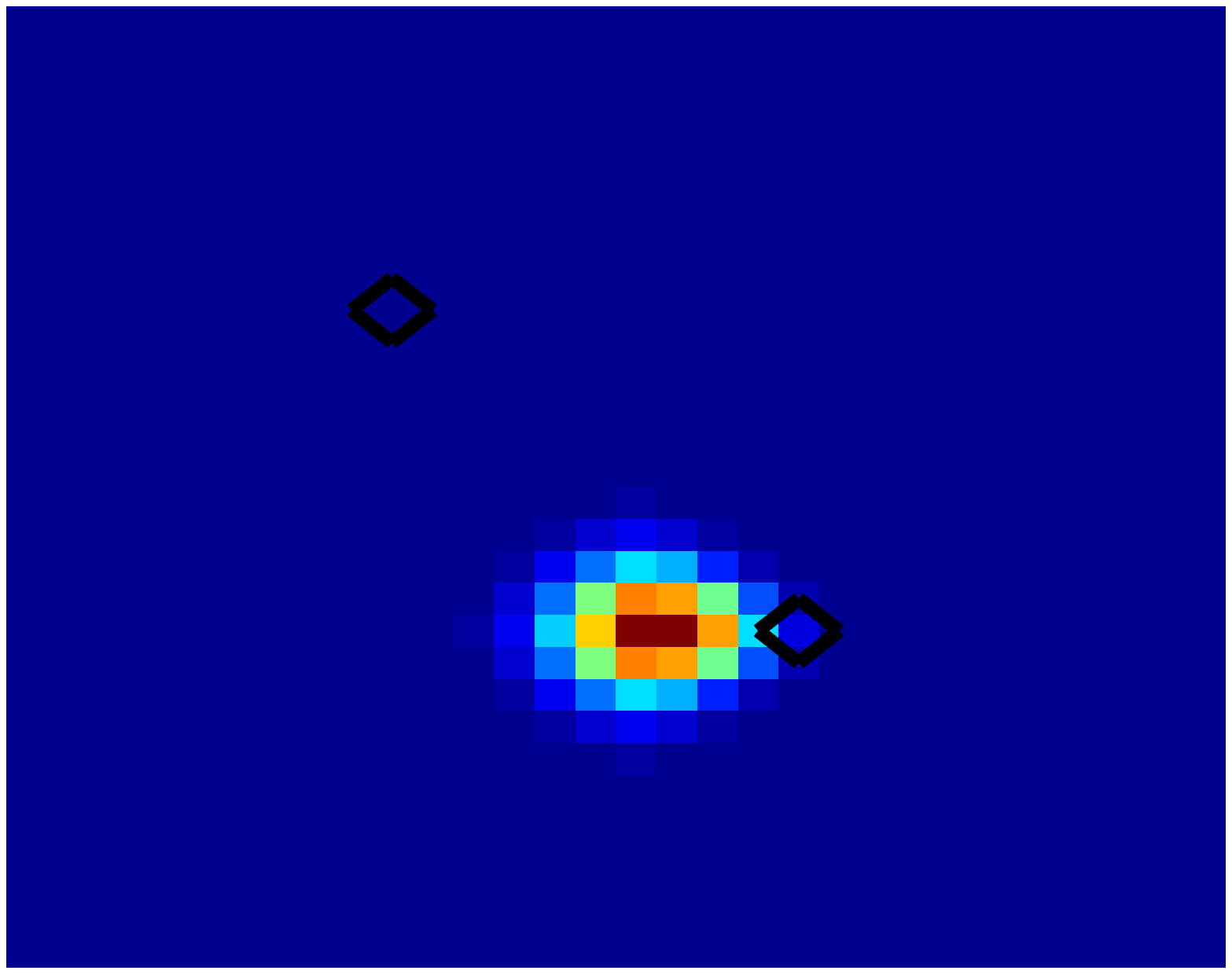}
    \includegraphics[width=1.0in]{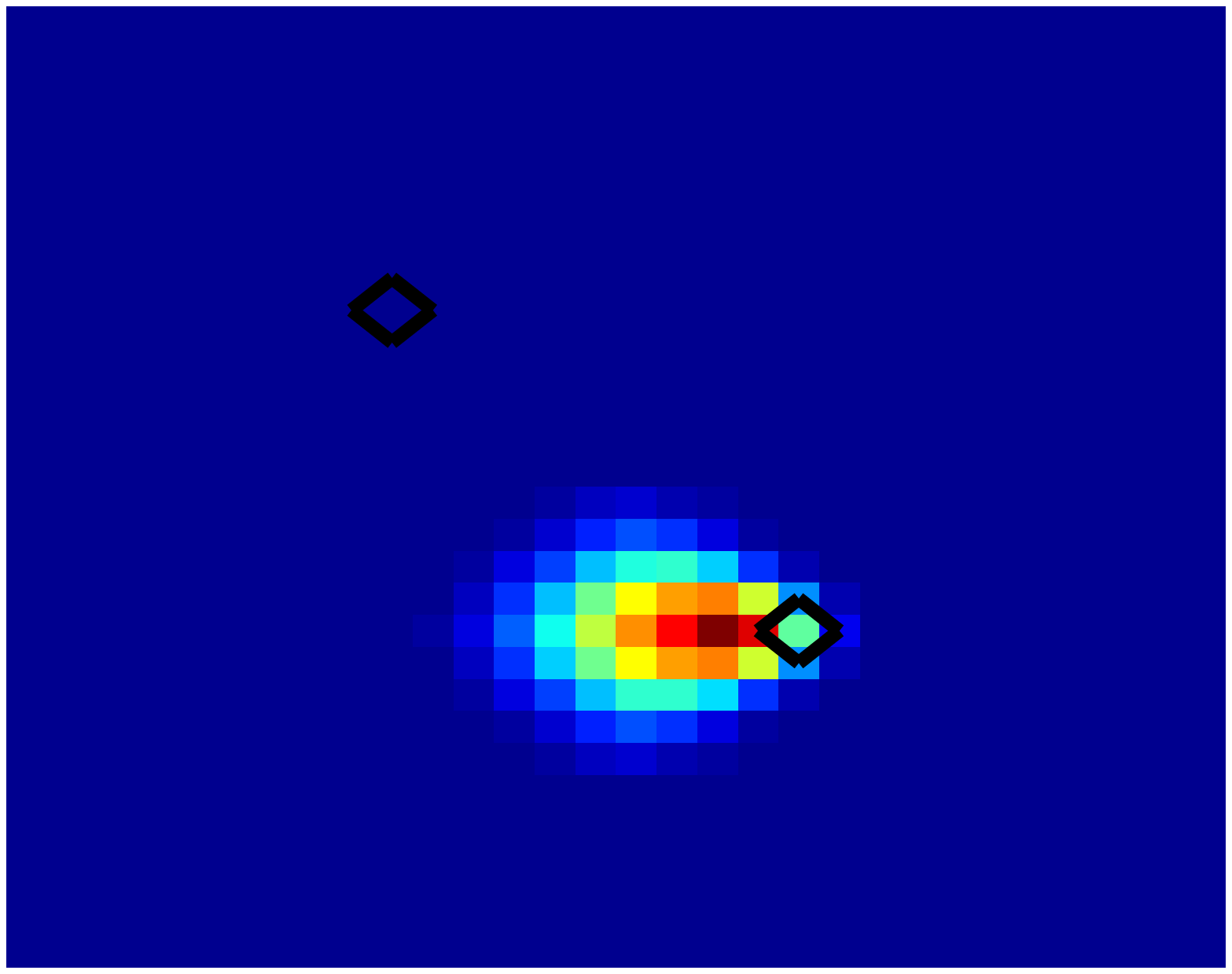}
    \includegraphics[width=1.0in]{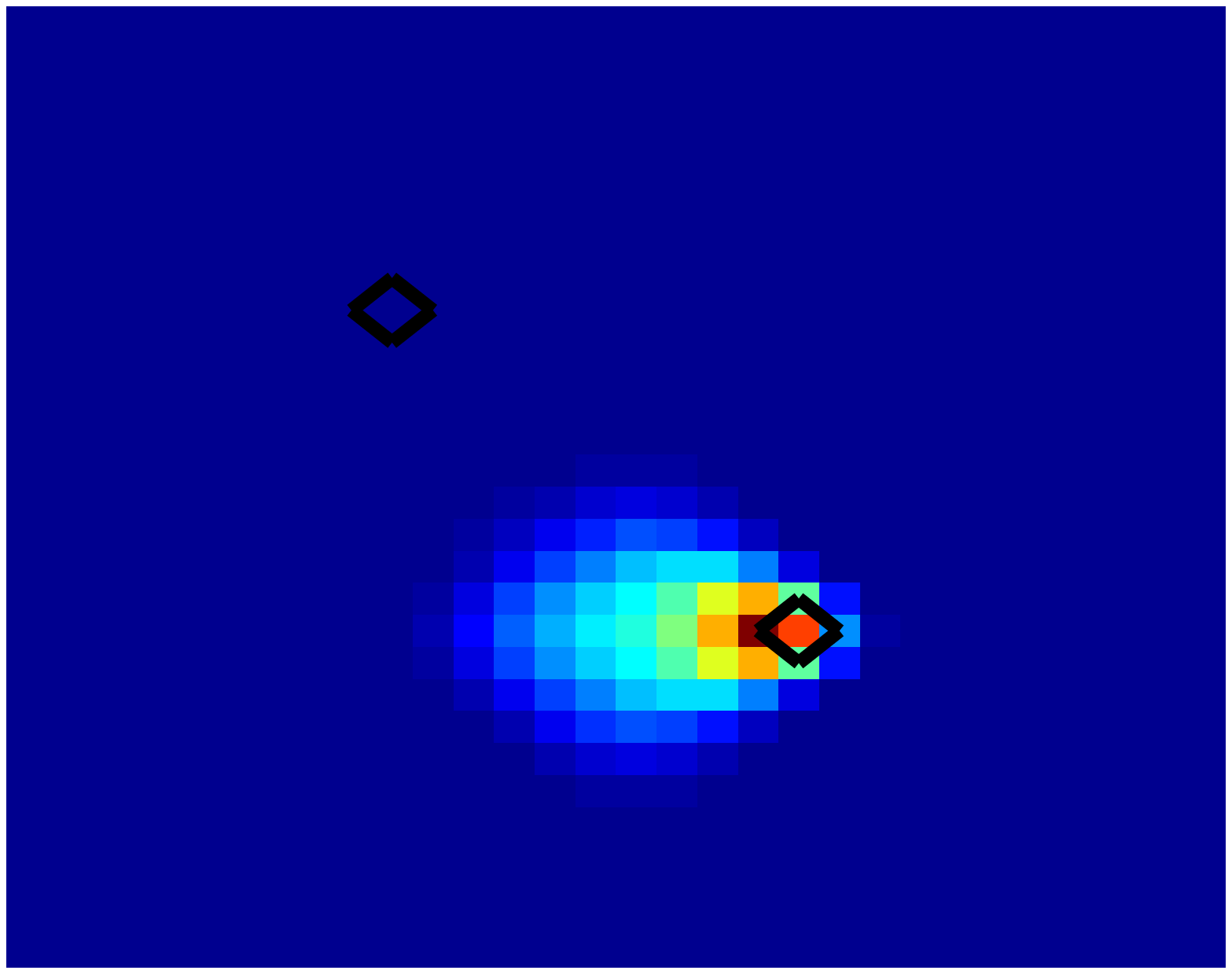}
    \includegraphics[width=1.0in]{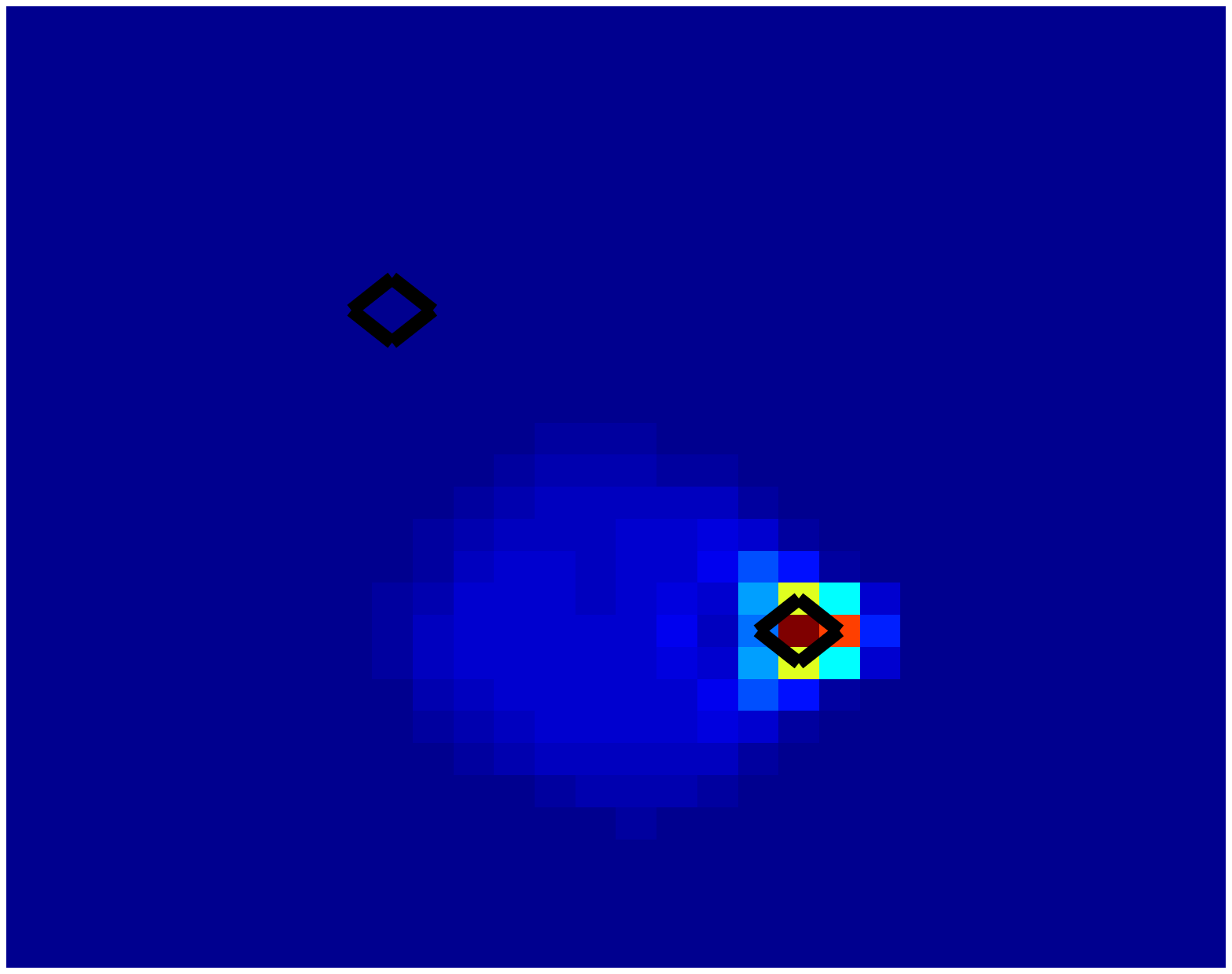}\\
    \includegraphics[width=1.0in]{Figs/Sparcify2_f.eps}\quad
    \includegraphics[width=1.0in]{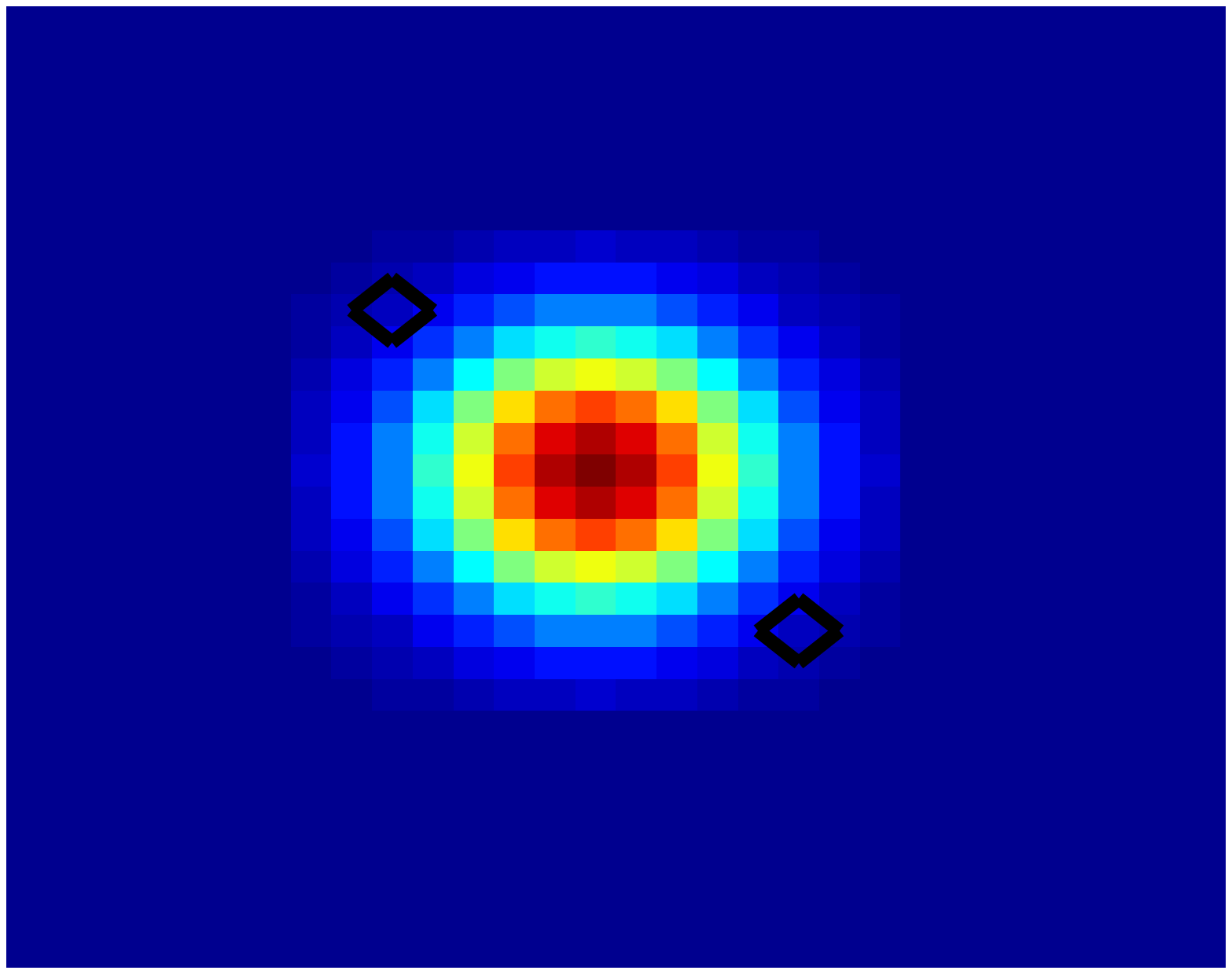}
    \includegraphics[width=1.0in]{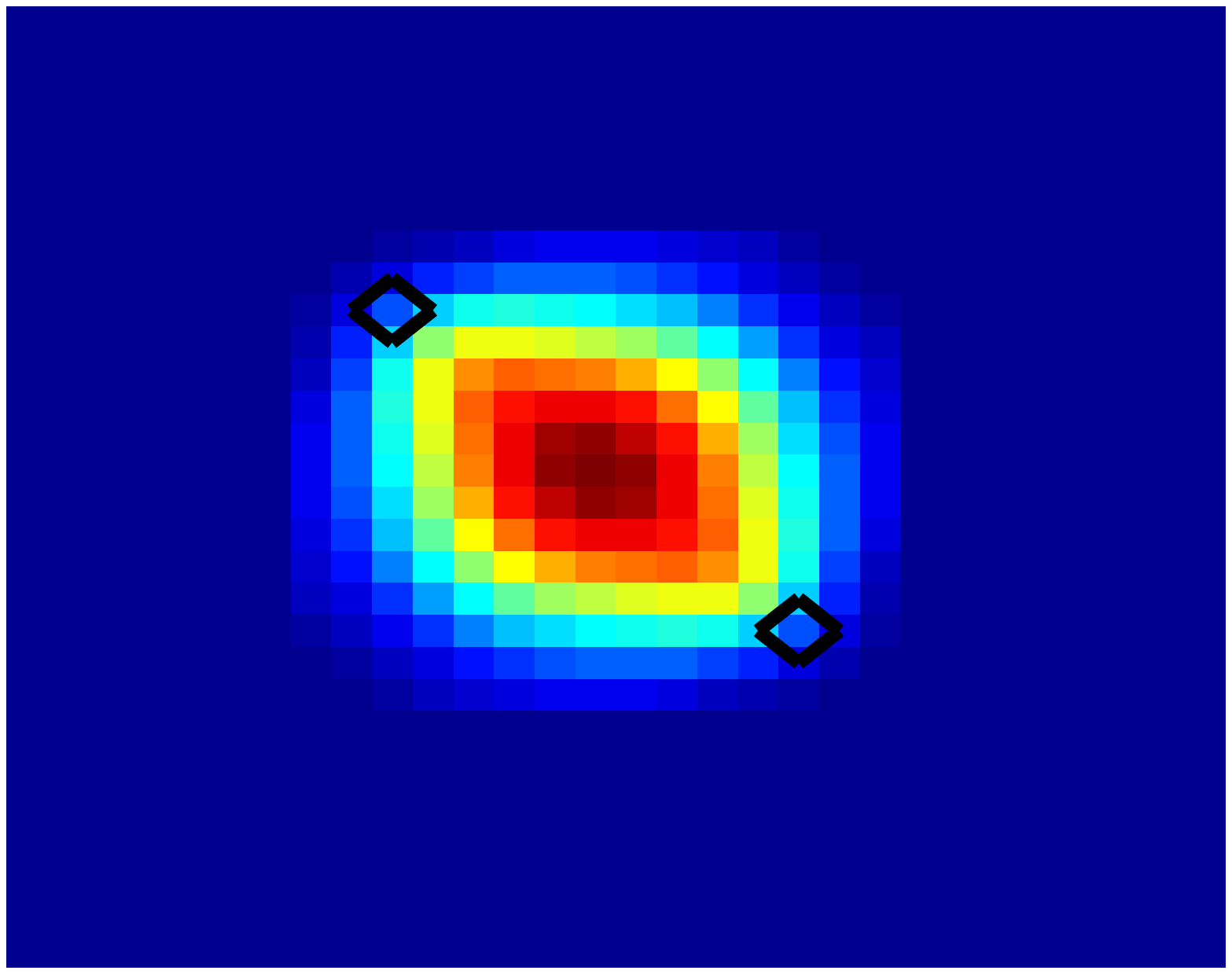}
    \includegraphics[width=1.0in]{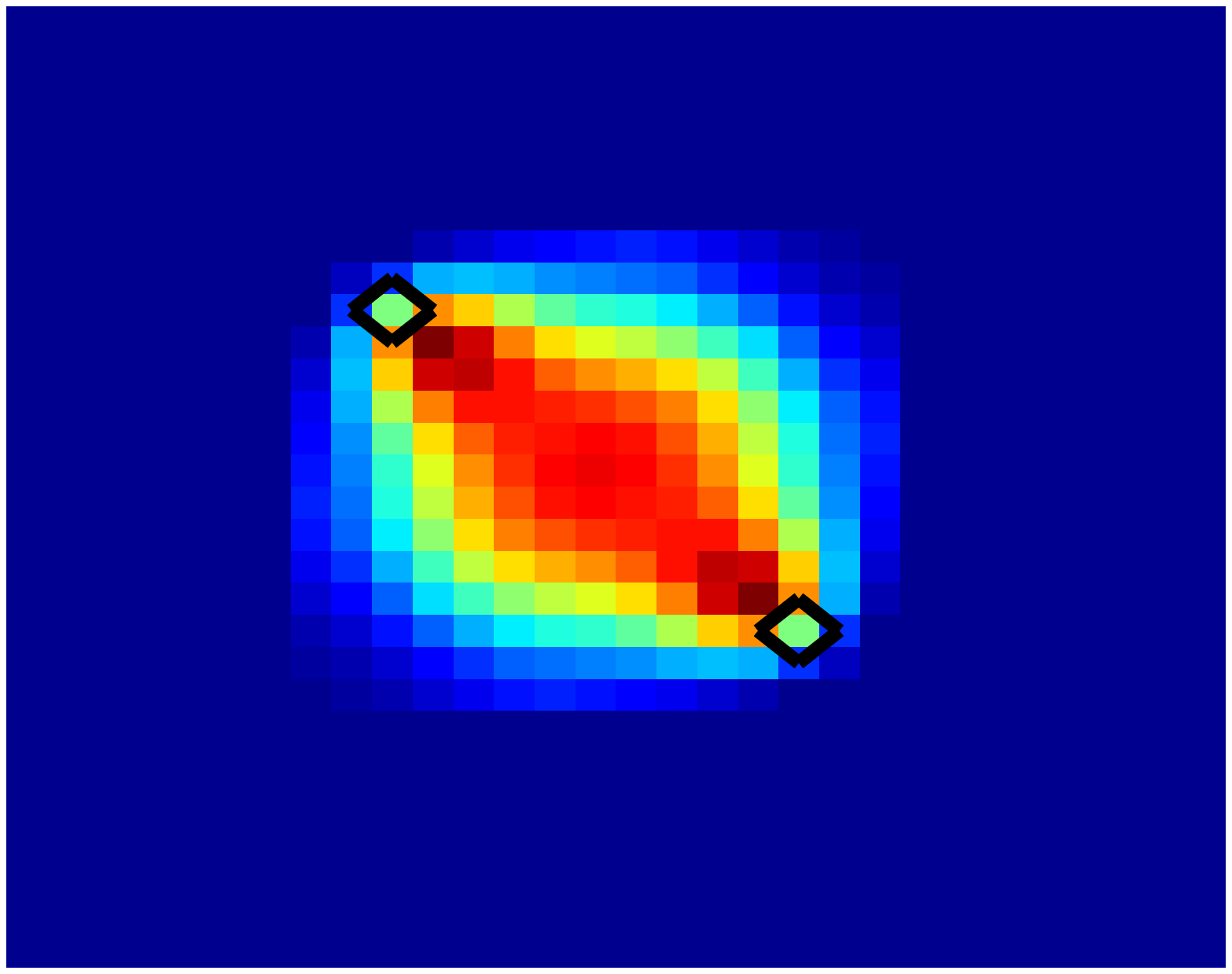}
    \includegraphics[width=1.0in]{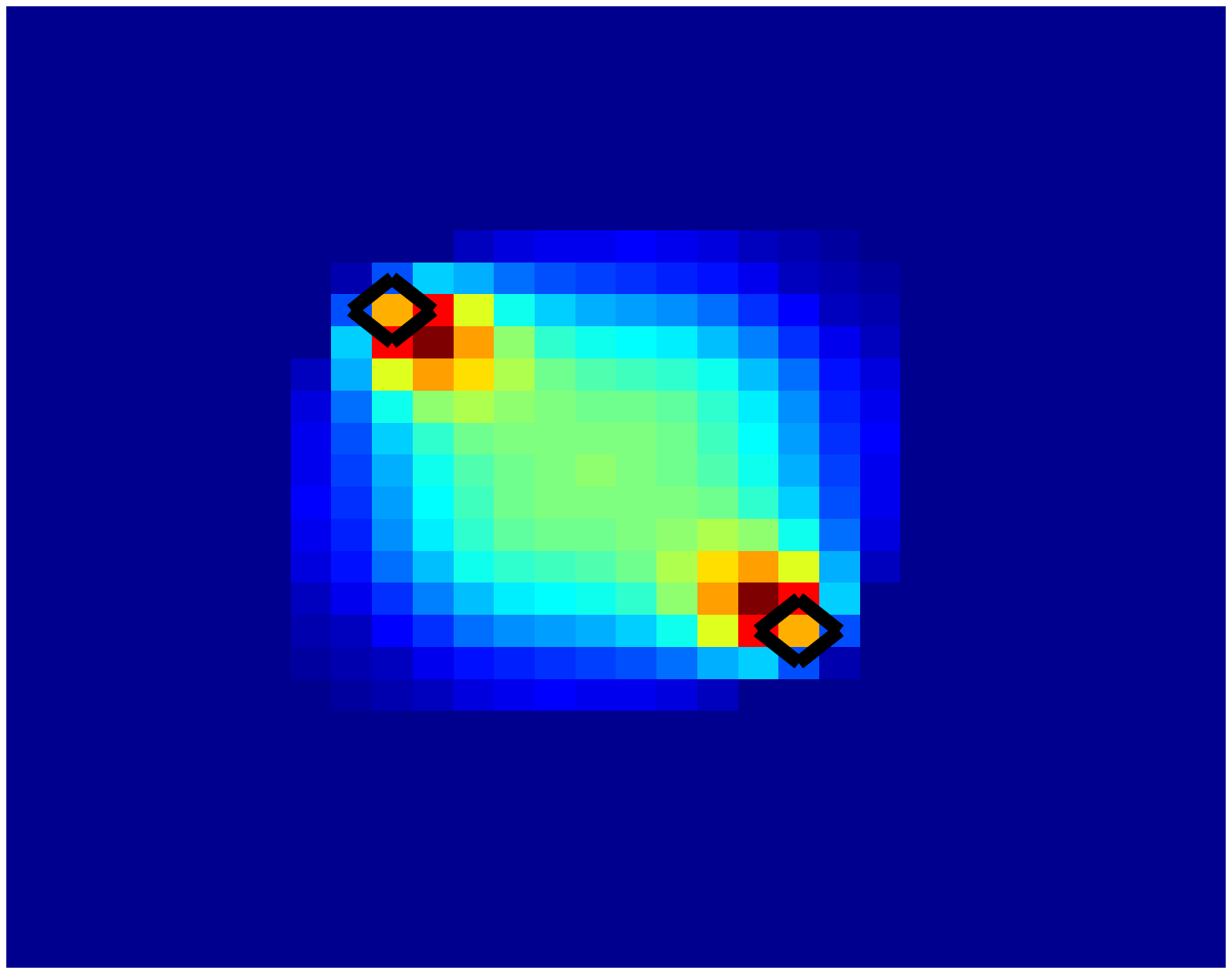}
    \includegraphics[width=1.0in]{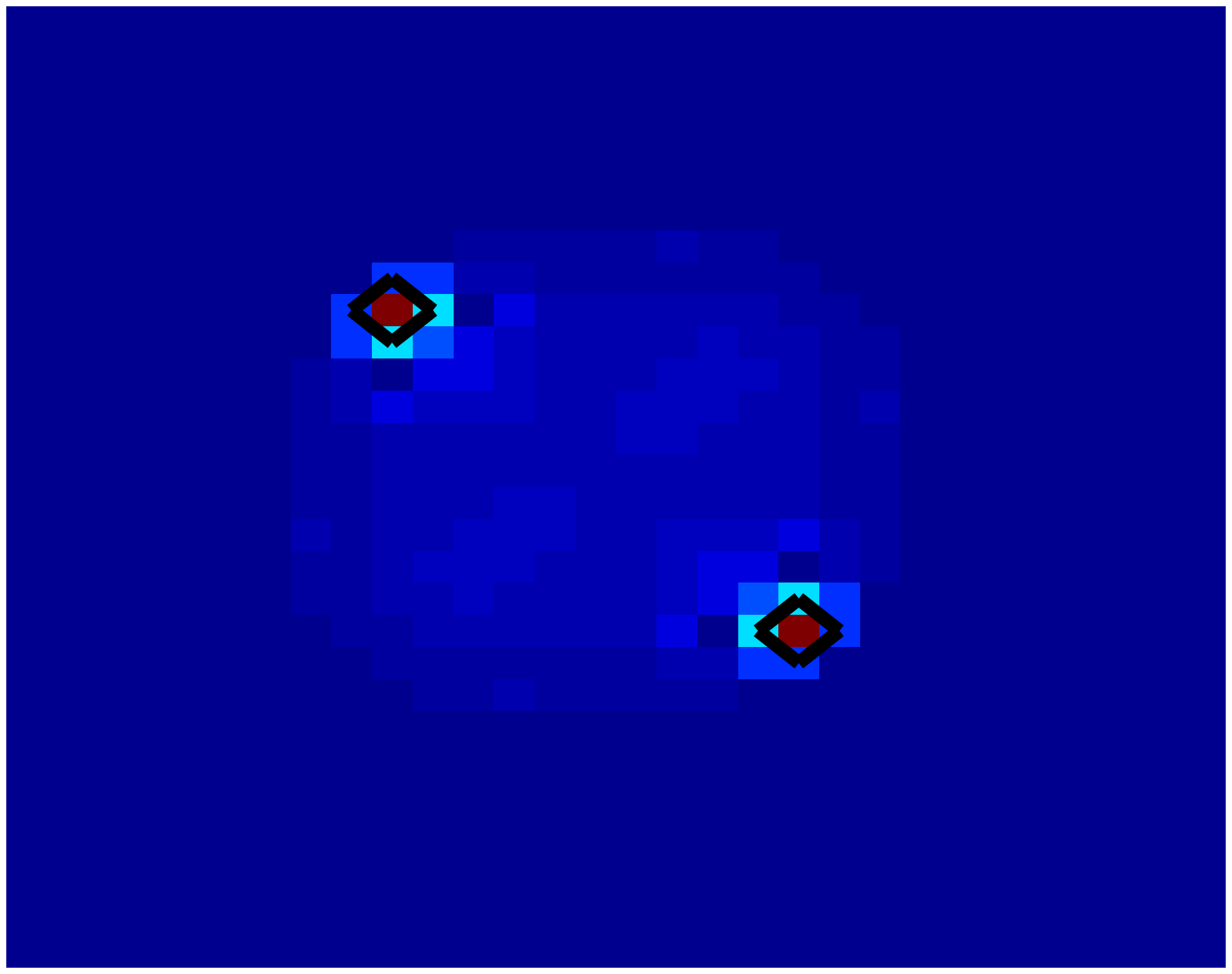}
    \caption{The first column are observed images $f$ obtained from convolution of a Gaussian with one or two spikes. The the second column shows the initial functions $u_0$ and the remaining 4 columns show the corresponding solutions of spatial motion PDE \eqref{motionPDE} at time T=250, 500, 750, 2000 respectively. The difference between the second and third rows reflect the fact that with the same $f$, different choices of initial value will lead to different results.}\label{Fig:Sparcify}
\end{figure}

Another and yet more important property of the motion PDE
\eqref{CauchyPDE} is that the solution $u(x,t)$ becomes sparser
asymptotically in $t$. We illustrate this phenomena numerically in
Figure \ref{Fig:Sparcify}, where the peaks of the solution $u$
move to the correct locations, and when they are close to the
desired locations, $u$ begins to become sparser which is desired
for sparse reconstruction problems. This explains why we can get
more accurate solutions when combining the motion PDE
\eqref{CauchyPDE} with FPC Bregman iterations (see Section
\ref{Results:accuracy} for more details).

It is also worth noticing that, comparing the second and the third
rows of Figure \ref{Fig:Sparcify}, for different initial functions
$u_0$, the solutions of PDE \eqref{CauchyPDE} have rather different
asymptotic behavior. This means that the residual $H(u)$ is not
convex under the Wasserstein metric. In fact, the solution in the
second row of Figure \ref{Fig:Sparcify} only approximates a local
minimizer of $H(u)$ with respect to the Wasserstein metric but not
a global one, i.e. $H(u)\neq 0$ in this case, but $\text{grad}_g H(u)=0$.

\section{Numerical Schemes and Experiments}\label{numerical}

\subsection{Numerical Implementation of PDE \eqref{CauchyPDE}}\label{implementation}

The finite difference scheme that we use for our motion
PDE is as follows:
\begin{equation}\label{Finite:Difference:Scheme}
\frac{u_{i,j}^{n+1}-u_{i,j}^{n}}{dt}=D^x_-\left[
\left|u_{i+\frac{1}{2},j}^n\right| D^x_+ p^n_{i,j} \right] +
D^y_-\left[ \left|u_{i,j+\frac{1}{2}}^n\right| D^y_+ p^n_{i,j}
\right],
\end{equation}
where \benn D^x_-u_{i,j}= u_{i,j}-u_{i-1,j},\quad D^x_+u_{i,j}=
u_{i+1,j}-u_{i,j}, \eenn \benn D^y_-u_{i,j}=
u_{i,j}-u_{i,j-1},\quad D^y_+u_{i,j}= u_{i,j+1}-u_{i,j} \eenn
\benn \left|u_{i+\frac{1}{2},j}^n\right|:=\left\{\begin{array}{ll}
\frac{|u^n_{i+1,j}|+|u^n_{i,j}|}{2}, &
u_{i+1,j}u_{i,j}>0,\\
0 & \mbox{otherwise}
\end{array}\right. \eenn
\benn \left|u_{i,j+\frac{1}{2}}^n\right|:=\left\{\begin{array}{ll}
\frac{|u^n_{i,j+1}|+|u^n_{i,j}|}{2}, &
u_{i,j+1}u_{i,j}>0,\\
0 & \mbox{otherwise}
\end{array}\right. \eenn
and
$$p^n=A^\top(Au^n-f)$$ where the operator $A$ is implemented by standard matrix
multiplication. The special definition of
$|u_{i+\frac{1}{2},j}^n|$ and $|u_{i,j+\frac{1}{2}}^n|$ above is
to make sure that the $\ell_1$ norm of $u^k$ is exactly preserved
in the discrete setting.

To guarantee the stability of the numerical solution, the time
step $dt$ should be chosen appropriately. Although it is hard to
derive a explicit stable condition for our nonlinear PDE, for each
step we can treat the PDE as a transport equation with a fixed
velocity field, and find a necessary stability condition, which is
given by

\be\label{CFL}
\frac{1}{dt^n}\geq \max_{i,j}\{|D_+^xp^n_{i,j}|,|D_+^yp^n_{i,j}|\}.
\ee

\subsection{Algorithm}\label{algorithm}

As we explained in section \ref{properties}, the PDE
\eqref{CauchyPDE} does not necessarily lead to a correct solution
of the optimization problem \eqref{L1:Minimization}. However, it
helps the function to be more sparse while decreasing the
residual. Therefore it can be combined with other existing
deconvolution algorithm and improve the performance. Since the
motion PDE \eqref{CauchyPDE} is time-dependent, it is natural to
combine  the numerical scheme \eqref{Finite:Difference:Scheme}
with FPC Bregman iteration \eqref{FPC+Bregman} and introduce the
following algorithm:
\begin{equation}\label{PDE+FPC+Bregman}
\begin{array}{ll}
 u^{k+1,N}& \leftarrow \left\{\begin{array}{ll} u^{k,l+\frac{1}{3}} & \leftarrow u^{k,l}-\delta A^\top(Au^{k,l}-f^k)\\
 u^{k,l+\frac{2}{3}} & \leftarrow u^{k,l+\frac{1}{3}}
+dt\nabla_x\cdot\left[ |u^{k,l+\frac{1}{3}}|\nabla_x\left(
A^\top(Au^{k,l}-f^k) \right)
\right]\qquad\text{(PDE\ Update)}\\
u^{k,l+1} & \leftarrow \text{shrink}(u^{k,l+\frac{2}{3}},\mu) \end{array}\right.\\
f^{k+1} &= f^k + f-Au^{k+1,N}.
\end{array}
\end{equation}
We shall refer to the above algorithm as PDE-FPC Bregman
iterations. 

We further note that for each iteration the PDE update only
slightly increases the complexity, because we do not need to
recalculate the term $A^\top(Au-f)$ and the finite difference
operator \eqref{Finite:Difference:Scheme} can be implemented
rather efficiently.

We want to point out that this PDE update step can be similarly
imposed to many other iterative thresholding based algorithms,
e.g. the linearized Bregman iterative method \cite{OMDY} or the
BOS method \cite{zhang2009bregmanized}. As far as we tested, this
complementary step always improves the performance and accuracy of
these algorithms without increasing the computational complexity
very much. In many cases the improvement is very significant.

\textbf{Remark:}
Note that the inner loop of \eqref{PDE+FPC+Bregman}
contains $N\ge1$ iterations, as in the original FPC algorithm. The choice of $N$ is not critical for the convergence of the scheme, however it affects the performance. When $N$ is small, the algorithm converges faster (in terms of total number of iterations) while the decay of the residual is more oscillatory; when $N$ is large, it converges slower while the decay of the residual appears more stable. In our experiments we choose $N=10$, which seems to give us a good balance between the speed and decay of residual.

\subsection{Numerical Results}\label{results}
\begin{figure}[htp]
\centering
    \includegraphics[width=2.0in]{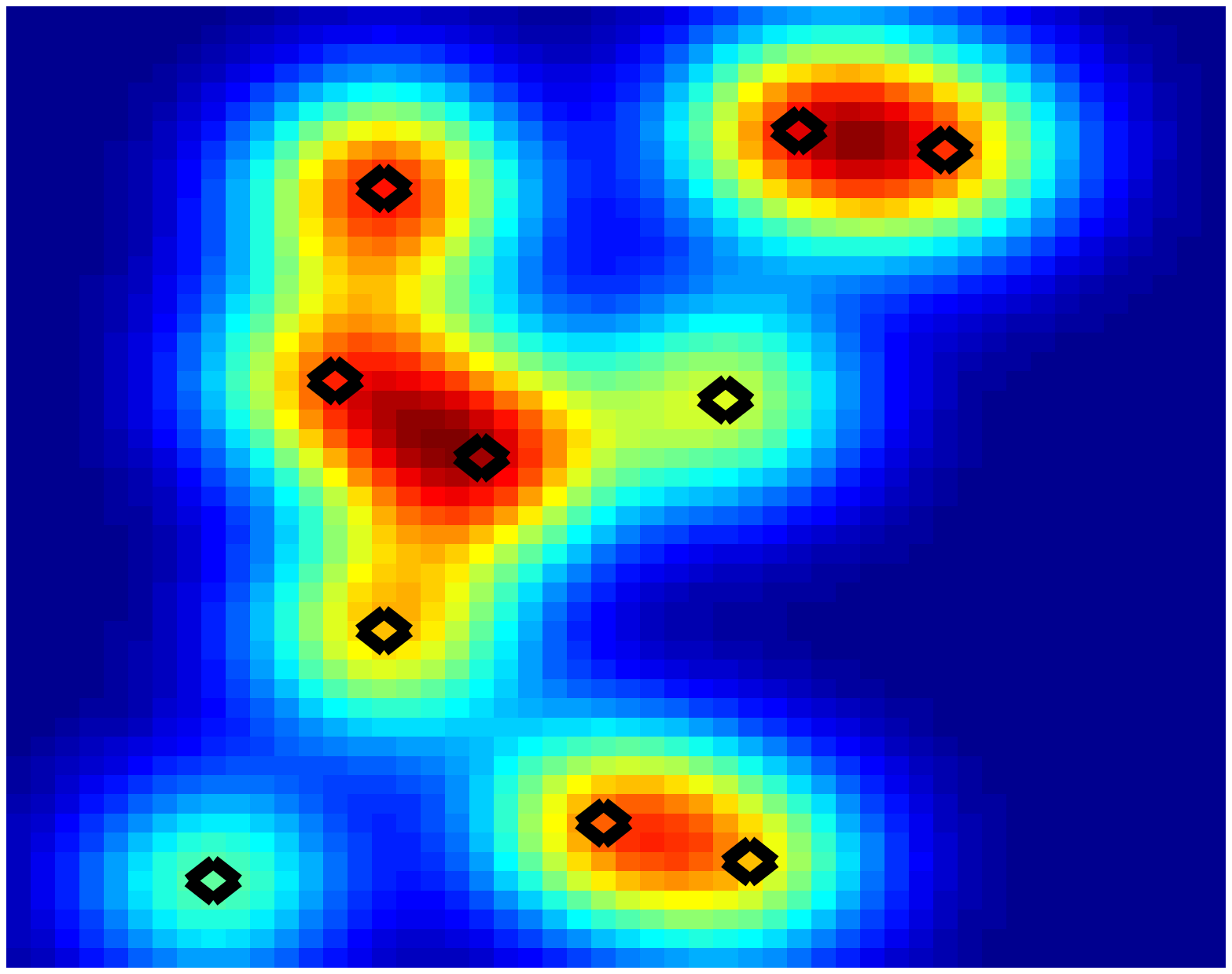}
    \includegraphics[width=2.0in]{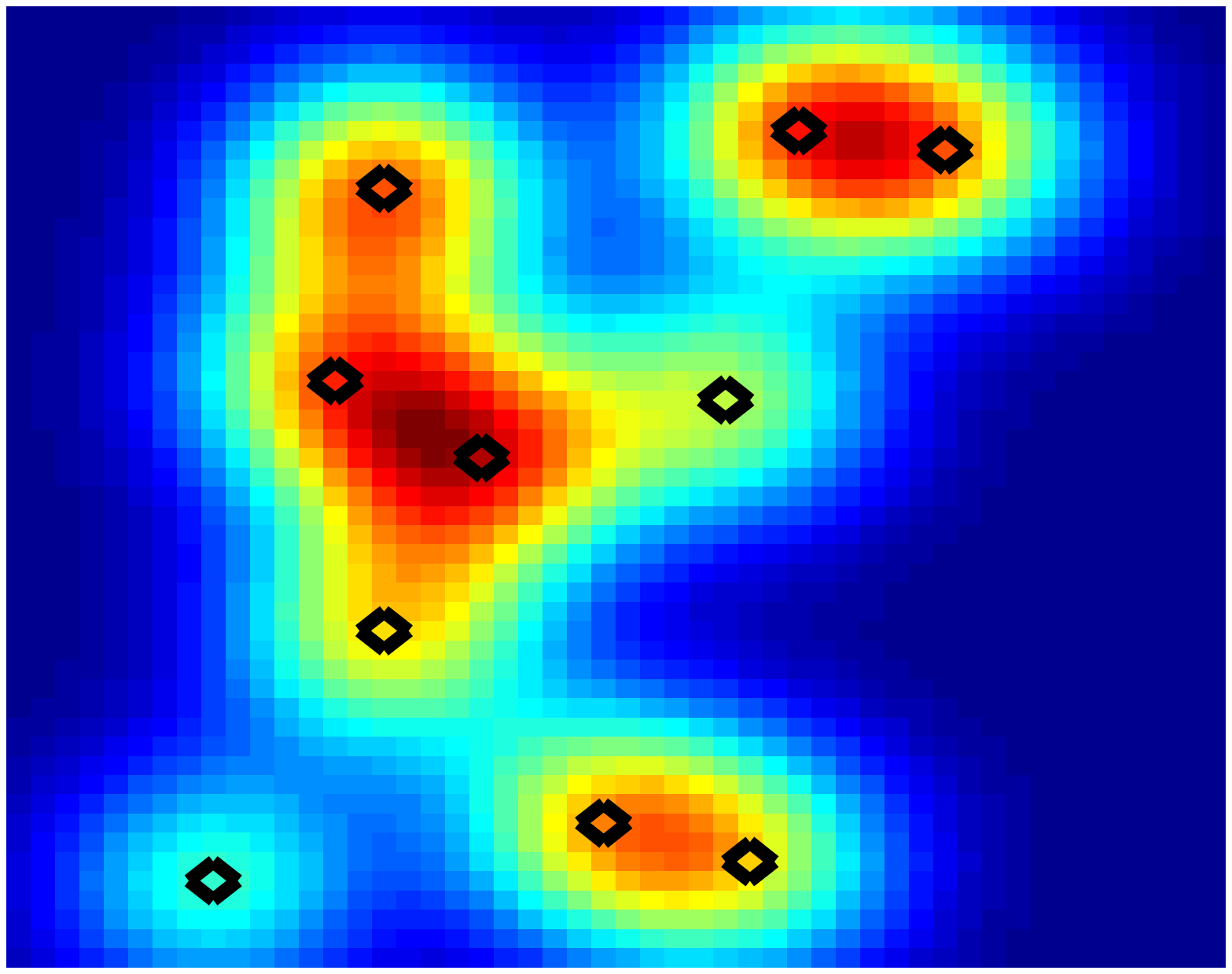}
    \includegraphics[width=2.0in]{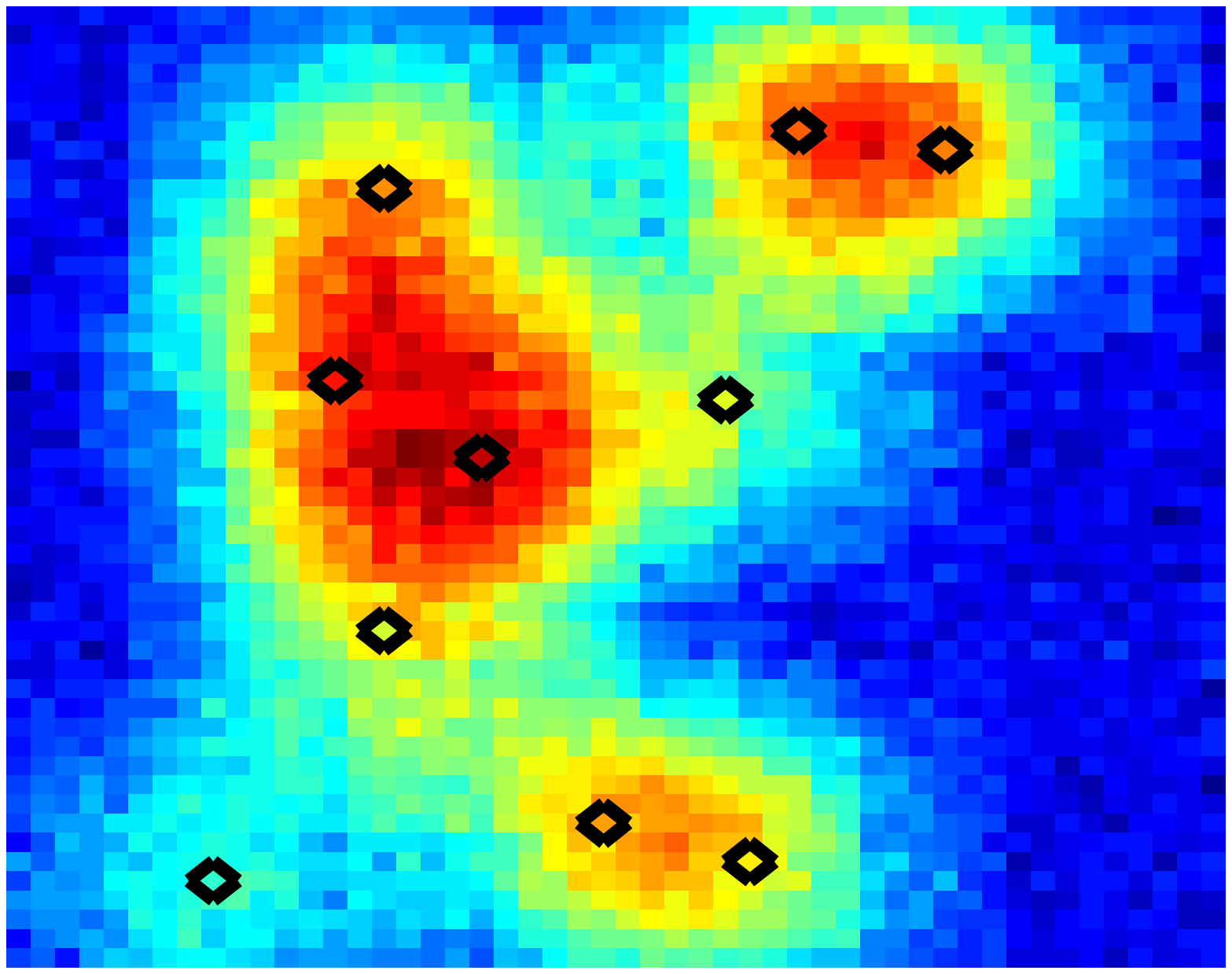}
    \caption{Clean observed images generated by different Gaussian kernels(left and middle). Noisy observed image (right) with SNR=15.87. Black circles indicating
    the locations of the spikes.}\label{Fig:CleanNoisy:f}
\end{figure}

In this section, we compare the performance of FPC Bregman
interations \eqref{FPC+Bregman} with that of PDE-FPC Bregman
iterations \eqref{PDE+FPC+Bregman}. We will first consider the
case where there is not any noise in $f$, and compare the
convergence speed. Then we consider a noisy case and compare the
accuracy of the solutions obtained from the two algorithms. The
operator $A$ is taken to be a convolution operator with a Gaussian
kernel. The clean and noisy $f$ are shown in Figure
\ref{Fig:CleanNoisy:f}, where the size of the images is
$50\times50$. The two kernels of the clean blurry images are
generated by Matlab function \texttt{fspecial('gaussian',[41
41],$\sigma$)} with $\sigma=4$ and 4.5 respectively. The noisy
images are generated by adding Gaussian white noise to the clean
blurry image with kernel \texttt{fspecial('gaussian',[21 21],$5.5$)}.
The intensities of the spikes are generated randomly in $[0.5,1]$
with locations indicated by the black circles (Figure
\ref{Fig:CleanNoisy:f}).

\subsubsection{Convergence Speed}\label{Results:speed}
We consider the noise free case (left two images in Figure
\ref{Fig:CleanNoisy:f}). Here we take the stopping criteria
$\frac{\|u^k-\bar u\|_2}{\|\bar u\|_2}<10^{-2}$, where $\bar u$ is
the true solution. The following table summaries the comparison
results with various choices of parameter $\mu$. The time step
$dt$ for the PDE is chosen to reflect best performance for given
$\mu$. From Table \ref{Table:Speed1} and Table \ref{Table:Speed2}
we can see that PDE-FPC Bregman is generally faster than FPC
Bregman in terms of computational time. Furthermore, the smaller
$\mu$ gets, the more PDE-FPC Bregman outperforms FPC Bregman. The
reason is that for small $\mu$, $u^k$ is more regular (i.e. less
sparse) than for large $\mu$ due to thresholding, and hence the
PDE update in \eqref{PDE+FPC+Bregman} will have a much nicer
behavior. Also, a larger time step $dt$ is allowed when $u^k$ is
more regular.

It is also worth noticing from Table \ref{Table:Speed1} and
\ref{Table:Speed2} that the best performance of PDE-FPC Bregman is
6 and 10 times faster respectively than that of FPC Bregman
(marked by ``$\spadesuit$"). The more blurry is the image, the
bigger is the improvement of PDE-FPC Bregman over FPC Bregman.

\begin{table}[ht]
\caption{Comparisons of FPC Bregman \eqref{FPC+Bregman} and
PDE-FPC Bregman \eqref{PDE+FPC+Bregman}. The symbol
``$\spadesuit$" labels the best results for FPC Bregman and
PDE-FPC Bregman.}\label{Table:Speed1}
\begin{center}
\hspace*{-0.25in}
\begin{tabular}{r|r||c|c||c|c}\hline\hline
 & & \multicolumn{2}{c||}{FPC Bregman} & \multicolumn{2}{c}{PDE-FPC Bregman} \\\hline
 $(\delta, \mu)$&$dt$ & \#Iterations & Time (sc) & \#Iterations & Time (sc) \\\hline
 (2, 0.5) & 0.2 & 2824 & 13.21 & 2816 & 14.97 \\\hline
 (2, 0.2) & 0.5 & 2196 & 10.69 $\spadesuit$ & 1705 & 9.23 \\\hline
 (2, 0.1) & 2.0 & 2585 & 12.43 & 1324 & 7.06 \\\hline
 (2, 0.05) & 4.5 & 3365 & 16.04 & 1291 & 6.91 \\\hline
 (2, 0.01) & 300.0 & 4997 & 24.20 & 587 & 3.20 \\\hline
 (2, 0.002) & 800.0 & 9622 & 47.76 & 280 & 1.64 \\\hline
 (2, 0.001) & 1200.0 & 13386 & 67.73 & 243 & 1.48 $\spadesuit$ \\\hline
\end{tabular}
\end{center}
\end{table}

\begin{table}[ht]
\caption{Comparisons of FPC Bregman \eqref{FPC+Bregman} and
PDE-FPC Bregman \eqref{PDE+FPC+Bregman}. The symbol
``$\spadesuit$" labels the best results for FPC Bregman and
PDE-FPC Bregman.}\label{Table:Speed2}
\begin{center}
\hspace*{-0.25in}
\begin{tabular}{r|r||c|c||c|c}\hline\hline
 & & \multicolumn{2}{c||}{FPC Bregman} & \multicolumn{2}{c}{PDE-FPC Bregman} \\\hline
 $(\delta, \mu)$&$dt$ & \#Iterations & Time (sc) & \#Iterations & Time (sc) \\\hline
 (2, 0.5) & 0.25 & 5327 & 25.37 & 4808 & 25.33 \\\hline
 (2, 0.2) & 0.9 & 4115 & 19.81 $\spadesuit$ & 2904 & 15.26 \\\hline
 (2, 0.1) & 2.0 & 4651 & 22.13 & 2562 & 13.70 \\\hline
 (2, 0.05) & 5.0 & 5836 & 28.31 & 1876 & 10.31 \\\hline
 (2, 0.01) & 300.0 & 8460 & 41.42 & 1030 & 5.61 \\\hline
 (2, 0.002) & 1500.0 & 14996 & 74.01 & 446 & 2.53 \\\hline
 (2, 0.001) & 2000.0 & 22466 & 112.29 & 328 & 1.91 $\spadesuit$ \\\hline
\end{tabular}
\end{center}
\end{table}

\subsubsection{Accuracy}\label{Results:accuracy}

We now consider the noisy case (right image in Figure
\ref{Fig:CleanNoisy:f}). Experimental results are summarized in
the following Table \ref{Table:Accuracy}. Since there is noise in
$f$, the error, which is taken to be $\frac{\|u^k-\bar
u\|_2}{\|\bar u\|_2}$, generally decreases first and then
increases, and $u^k$ will eventually converge to something noisy.
Therefore, in order to reflect the best performance of both
algorithms, we recorded the number of iterations and computation
time that correspond to the smallest error possible for each given
set of parameters. From Table \ref{Table:Accuracy} we can see
that, for each given $\mu$, although PDE+FPC Bregman may not
always be faster than FPC Bregman, it always produces a solution
with a smaller error. If we consider the best set of parameters
for the two algorithms, PDF-FPC Bregman outperforms FPC Bregman in
terms of both computation time and accuracy.

\begin{table}[ht]
\caption{Comparisons of FPC Bregman \eqref{FPC+Bregman} and
PDE-FPC Bregman \eqref{PDE+FPC+Bregman}. The symbol
``$\spadesuit$" labels the best computational time and error
respectively for FPC Bregman and PDE-FPC
Bregman.}\label{Table:Accuracy}
\begin{center}
\hspace*{-0.25in}
\begin{tabular}{r|r||c|c|c||c|c|c}\hline\hline
 & & \multicolumn{3}{c||}{FPC Bregman} & \multicolumn{3}{c}{PDE-FPC Bregman} \\\hline
 $(\delta, \mu)$&$dt$ & \#Iterations & Time (sc) & Error & \#Iterations & Time (sc) & Error \\\hline
 (2, 0.8) & 0.01 & 13299 & 17.52 & 1.12e-001 & 13417 & 24.61 & 1.04e-001 \\\hline
 (2, 0.5) & 0.05 & 6086 & 7.74 & 2.68e-002 $\spadesuit$ & 6033 & 10.81 & 1.63e-002 \\\hline
 (2, 0.2) & 0.5 & 3663 & 4.61 & 5.70e-002 & 3321 & 5.94 & 2.40e-002 \\\hline
 (2, 0.05) & 1.5 & 2127 & 2.94 $\spadesuit$ & 4.87e-002 & 1424 & 2.59 & 1.42e-002 \\\hline
 (2, 0.01) & 3.0 & 2404 & 3.49 & 5.17e-002 & 1651 & 3.09 & 1.47e-002 \\\hline
 (2, 0.005) & 4.0 & 2511 & 3.81 & 9.00e-002 & 1406 & 2.73 & 1.93e-002 \\\hline
 (2, 0.002) & 45.0 & 2985 & 4.92 & 1.87e-001 & 571 & 1.20 $\spadesuit$ & 1.26e-002 $\spadesuit$ \\\hline
 (2, 0.001) & 40.0 & 2730 & 4.70 & 4.05e-001 & 623 & 1.38 & 5.22e-002 \\\hline
\end{tabular}
\end{center}
\end{table}

\section{Conclusions}\label{conclusion}

In this paper, we introduced a new nonlinear evolution PDE for
sparse deconvolution problems. We showed that this spatial motion
PDE preserves the $\ell_1$ norm while lowering the residual
$\|Au-f\|_2$. We also showed numerically that the solution of the
motion PDE tends to become sparse and converges to delta functions
asymptotically. Therefore, utilizing these properties of the PDE,
our proposed PDE-FPC Bregman iterations \eqref{PDE+FPC+Bregman}
outperformed the original FPC Bregman iterations
\cite{Yin:2008p611} in terms of both convergence speed and
reconstruction quality. This was strongly supported by our
numerical experiments. We finally note that our PDE can also be
used to enhance performance for many other iterative thresholding
based algorithms, e.g.
\cite{OMDY,zhang2009bregmanized,figueiredo2003algorithm}.

\bibliography{biblist}
\end{document}